\documentclass[preprint,11pt]{elsarticle}

\usepackage{amsmath,amsthm,amsfonts,amssymb,bm,
graphicx,algorithm,algpseudocode,hyperref,cleveref,xcolor,
subcaption,booktabs}

\crefname{algorithm}{Algorithm}{Algorithms} 
\Crefname{algorithm}{Algorithm}{Algorithms} 

\numberwithin{equation}{section}

\crefname{equation}{}{}
\crefname{algorithm}{Algorithm}{}
\crefname{theorem}{Theorem}{}
\crefname{remark}{Remark}{}
\crefname{figure}{Figure}{}

\begin{document}

\title{A time-splitting Fourier pseudospectral method for the Wigner(-Poisson)-Fokker-Planck equations}

\author[1,2]{Qian Yi}
\ead{yiqianqianyi@163.com}

\author[3]{Limin Xu\corref{cor1}}
\ead{xulimin@westlake.edu.cn}

\cortext[cor1]{Corresponding author}

\address[1]{ Department of Mathematics, Shantou University, Shantou, 515821, China}
\address[2]{ School of Mathematics and Statistics, Guangxi Normal University, Guilin, 541006, China}

\address[3]{Institute for Theoretical Sciences, Westlake Institute for Advanced Study, Westlake University, Hangzhou, 310030, China}

\begin{abstract}
In this article, we propose an efficient time-splitting Fourier pseudospectral method for the Wigner(-Poisson)-Fokker-Planck equations. The method achieves second-order accuracy in time and spectral accuracy in phase space, both of which are rigorously verified by numerical experiments. The validated scheme is then employed to study the long-time dynamics of these systems. We investigate the existence of steady states for both the Wigner-Fokker-Planck and Wigner-Poisson-Fokker-Planck equations. Notably, for the Wigner-Fokker-Planck system, our results provide numerical evidence for the existence of a steady state even when the external potential is far from harmonic. This is an important discovery, since this phenomenon has not been thoroughly established in theory.
\end{abstract}
\maketitle

\section{Introduction}

Quantum transport phenomena are pivotal in modern condensed matter physics, nanoscience, and semiconductor engineering \cite{datta1997electronic,divincenzo2000physical}. An accurate description of non-equilibrium quantum systems is essential for both understanding the microscopic properties of materials and developing novel quantum technologies. The Wigner function provides a unique and powerful framework for this purpose, offering a phase-space representation of quantum mechanics that connects directly to the language of classical kinetic theory \cite{Wigner1932quantum,hillery1984distribution}.

In realistic physical settings, quantum systems are rarely isolated. Interactions with an environment, such as a phonon bath or impurities, introduce dissipation and decoherence \cite{zurek2003decoherence,jacoboni2012monte}. To model such open quantum systems, the Wigner equation is coupled with a Fokker-Planck operator, which accounts for these environmental effects, resulting in the Wigner-Fokker-Planck (WFP) equation \cite{kolobov2003fokker}.
Furthermore, in systems with high particle density, such as in semiconductor devices, the electrostatic potential is significantly influenced by the charge distribution itself. Therefore, a self-consistent treatment is necessary, which is achieved by coupling the WFP equation with the Poisson equation \cite{degond1990quantum,zweifel1993poisson}. This leads to the more comprehensive Wigner-Poisson-Fokker-Planck (WPFP) equation, which provides a complete model for interacting, open quantum systems \cite{arnold2008mathematical}.

At present, the mathematical theory for the WFP and WPFP equations is well-established. For the WFP system, the long-time behavior of solutions has been analyzed \cite{sparber2004long}. For the more complex WPFP system, the existence, uniqueness, and asymptotic behavior of solutions have been rigorously proven \cite{arnold2004analysis,arnold2007dispersive,arnold2002periodic}, including for cases with singular potentials \cite{li2017wigner}.

Despite this extensive body of theoretical research, the development of numerical methods specifically for full WFP and WPFP systems remains relatively limited. Notable examples include the continued-fraction method developed by Palacios for the WFP equation \cite{garcia2004solving} and a discontinuous Galerkin scheme proposed by Gamba et al. \cite{gamba2009adaptable}.

The literature on numerical methods is far more developed for simplified versions of WFP/WPFP models. In the collisionless limit (i.e., without the Fokker-Planck operator), the equation becomes the Wigner(-Poisson) equation. A wide array of deterministic methods has been developed for this system, ranging from early finite difference schemes \cite{frensley1987wigner,frensley1990boundary,jensen1991methodology} to higher-order techniques such as spectral element methods \cite{shao2011adaptive}, WENO schemes \cite{dorda2015weno}, hyperbolic moment method \cite{li2014numerical} and operator-splitting spectral methods \cite{chen2019high,chen2022higher}. Alongside these, various stochastic methods have also been extensively explored\cite{shao2016computable,wagner2016random,sellier2014benchmark,rjasanow1996stochastic,muscato2010numerical,yan2015monte,shao2015comparison,muscato2016class,xiong2017wigner}.
Similarly, in the classical limit, the equation reduces to the Vlasov(-Poisson)-Fokker-Planck equation. The numerical literature for this system is also extensive, encompassing both deterministic particle methods \cite{wollman2005numerical,wollman1996numerical,havlak1998deterministic,wollman2009numerical}, stochastic particle methods \cite{allen1994computational,havlak1996numerical}, and grid-based methods, such as finite element \cite{asadzadeh2007convergence} and finite difference schemes \cite{schaeffer1998convergence}.

To address the aforementioned gap in the literature for the full WFP and WPFP systems, we propose an efficient and accurate time-splitting Fourier spectral (TSSP) method in this paper. Our approach employs operator splitting to decompose the governing equation into a set of simpler subproblems, each of which is then solved with a high precision spectral Galerkin or spectral collocation method. The resulting scheme is not only highly accurate, leveraging the spectral convergence of the spatial discretization and second-order accuracy of the splitting, but also computationally efficient and straightforward to implement.

The remainder of this paper is organized as follows. Section \ref{sec: Review of WPFP} provides a detailed overview of WFP and WPFP equation and its properties. Section \ref{sec:numerical_method} elaborates on the proposed TSSP method, detailing the discretization strategies for each subproblem. Section \ref{sec: numerical results} presents a series of numerical experiments to validate the accuracy and demonstrate the capabilities of our method. Finally, Section \ref{sec: conclusion} concludes the paper and discusses potential directions for future research.
\section{The Wigner(-Poisson)-Fokker-Planck equation}
\label{sec: Review of WPFP}

The dynamics of an open quantum system are described by the Wigner distribution function $W(\bm{x},\bm{\xi},t)$, whose evolution in the $2d$-dimensional phase space is governed by the dimensionless Wigner(-Poisson)-Fokker-Planck equation:
\begin{equation}
\begin{aligned}
    \frac{\partial W(\bm{x},\bm{\xi},t)}{\partial t}
    +\bm{\xi}\cdot \nabla_{\bm{x}}W+\Theta[V]W
    &=Q(W),\quad \bm{x},\bm{\xi}\in\mathbb{R}^d,\\
    W(\bm{x},\bm{\xi},0)&=W_0(\bm{x},\bm{\xi}),
\end{aligned}
\end{equation}
where $\bm{x} \in \mathbb{R}^d$ and $\bm{\xi} \in \mathbb{R}^d$ are the position and momentum vectors, respectively. The Wigner function $W(\bm{x},\bm{\xi},t)$ is normalized such that:
$$\iint_{\mathbb{R}^{d}\times \mathbb{R}^{d}} W(\bm{x},\bm{\xi},t)\mathrm{d}\bm{x}\mathrm{d}\bm{\xi}=1.$$

A key feature of the Wigner function is that it can attain negative values, and is thus considered a quasi-probability distribution.

The term $\Theta[V]W$ describes the interaction with the potential $V(\bm{x},t)$. The operator $\Theta[V]$ is a pseudo-differential operator, reflecting the inherently nonlocal nature of quantum mechanics, and is defined as:
    $$
    \begin{aligned}
        (\Theta[V]W)(\bm{x},\bm{\xi},t)&=\frac{1}{(\pi\varepsilon)^d}\int_{\mathbb{R}^d\times\mathbb{R}^d}\delta V(\bm{x},\bm{y},t)W(\bm{x},\bm{\xi}^{\prime},t)\mathrm{e}^{\frac{\mathrm{i}}{\varepsilon}2\bm{y}\cdot\left(\bm{\xi}^{\prime}-\bm{\xi}\right)}\mathrm{d}\bm{\xi}^{\prime}\mathrm{d}\bm{y},\\
        \delta V(\bm{x},\bm{y},t)&=\frac{\operatorname{i}}{\varepsilon}\left(V(\bm{x}+\bm{y},t)-V(\bm{x}-\bm{y},t)\right),
    \end{aligned}
    $$
    where $\varepsilon$ is the semi-classical parameter. The nature of the potential $V(\bm{x},t)$ distinguishes two important cases. When $V(\bm{x},t)$ is a prescribed external field, the equation is known as the Wigner-Fokker-Planck equation. When $V(\bm{x},t)$ is determined self-consistently through the particle density $\rho(\bm{x},t)$ via the Poisson equation, it is called as the Wigner-Poisson-Fokker-Planck  equation:
    $$
    \Delta V(\bm{x},t)=\alpha \rho(\bm{x},t),\quad \text{where} \quad \rho(\bm{x},t)=\int_{\mathbb{R}^d}W(\bm{x},\bm{\xi},t)\mathrm{d}\bm{\xi}.
    $$
    Here, $\rho(\bm{x},t)$ is the particle density at $\bm{x}$, and $\alpha=\pm 1$ indicates repulsive or attractive forces. In the classical limit, as $\varepsilon\to0$, the operator $\Theta[V]$ converges to the classical Liouville operator:
    \begin{align}\label{eq: nonlocal reduce to local}
     \Theta[V] W \rightarrow-\nabla_{\boldsymbol{x}} V(\boldsymbol{x}) \cdot \nabla_{\boldsymbol{\xi}} W.
    \end{align}
    
   The collision operator $Q(W)$ accounts for the dissipative and diffusive effects arising from the system's interaction with an external environment. Its general form is:
    $$
    Q(W)=D_{pp}\Delta_{\bm{\xi}}W+2\gamma \operatorname{div}_{\bm{\xi}}(\bm{\xi}W)+D_{qq}\Delta_{\bm{x}}W+2D_{pq}\operatorname{div}_{\bm{x}}(\nabla_{\bm{\xi}}W),
    $$
    where $\gamma\ge0$ is the friction parameter, and the non-negative coefficients $D_{pp}$, $D_{pq}$, and $D_{qq}$ form the phase-space diffusion matrix. These coefficients are model-dependent; a common example derived from a Markovian approximation of an electron in an oscillator bath yields:
    $$
    \gamma=\frac{\eta}{2},\quad D_{pp}=\frac{\eta}{\beta},\quad D_{qq}=\frac{\beta\eta\varepsilon^2}{12},\quad D_{pq}=\frac{\beta\Omega\eta\varepsilon^2}{12\pi},
    $$
    where $\eta$ is the bath coupling constant, $\Omega$ is the cutoff frequency, and $\beta=1/(k_B T)$ is the inverse temperature\cite{arnold2004analysis}. The terms involving $D_{qq}$ and $D_{pq}$ are often referred to as ``quantum diffusion" as they vanish in the classical limit ($\varepsilon \to 0$). In this limit, combined with \cref{eq: nonlocal reduce to local}, the WPFP equation reduces to the classical Vlasov-Poisson-Fokker-Planck equation.


The macroscopic behavior of the system can be characterized by several physical quantities, derived by taking moments of the Wigner function. The primary local quantities include the particle density $\rho(\bm{x},t)$, the current density $\bm{j}(\bm{x},t)$, and the energy density $e(\bm{x},t)$:
\begin{align*}
    \rho(\bm{x},t) &= \int W(\bm{x},\bm{\xi},t)\mathrm{d}\bm{\xi}, \\
    \bm{j}(\bm{x},t) &= \int \bm{\xi} W(\bm{x},\bm{\xi},t)\mathrm{d}\bm{\xi}, \\
    e(\bm{x},t) &= \int \left(\frac{|\bm{\xi}|^2}{2}+\widetilde{\alpha}V(\bm{x},t)\right)W(\bm{x},\bm{\xi},t)\mathrm{d}\bm{\xi},
\end{align*}
In the expression for energy density, the coefficient  $\widetilde{\alpha}$ depends on the nature of the potential $V(\bm{x},t)$. Specifically, $\widetilde{\alpha}=1$ for an external potential, whereas $\widetilde{\alpha}=\frac{1}{2}$  for a self-consistent potential.

Integrating these densities over the spatial domain yields the corresponding global quantities: the total number of particles $N(t)$, total momentum $\bm{J}(t)$, and total energy $E(t)$:
\begin{align*}
    N(t) &= \int \rho(\bm{x},t)\mathrm{d}\bm{x} = \iint W(\bm{x},\bm{\xi},t)\mathrm{d}\bm{x}\mathrm{d}\bm{\xi}, \\
    \bm{J}(t) &= \int \bm{j}(\bm{x},t)\mathrm{d}\bm{x} = \iint \bm{\xi} W(\bm{x},\bm{\xi},t)\mathrm{d}\bm{x}\mathrm{d}\bm{\xi}, \\
    E(t) &= \int e(\bm{x},t)\mathrm{d}\bm{x} = \iint \left(\frac{|\bm{\xi}|^2}{2}+\widetilde{\alpha}V(\bm{x},t)\right)W(\bm{x},\bm{\xi},t)\mathrm{d}\bm{x}\mathrm{d}\bm{\xi}.
\end{align*}

These macroscopic quantities are essential for characterizing the system's physical behavior. Their evolution reveals a fundamental distinction from closed quantum systems: for the WFP equation, while the total particle number is strictly conserved, the total energy are generally not, due to dissipative interactions with an environment. This is in stark contrast to ideal closed systems, where total energy is also a conserved quantity.

\section{Time-splitting Fourier pseudospectral method for the 1D WPFP equation}
\label{sec:numerical_method} 

In this section, we present the numerical scheme for solving the one-dimensional (1D) dimensionless WPFP equation. In 1D, the collision operator defined in \cref{sec: Review of WPFP} becomes:
\begin{equation}
	\begin{aligned}
		Q(W)=D_{pp}\frac{\partial^2 W}{\partial \xi^2}+2\gamma\frac{\partial (\xi W)}{\partial \xi}+D_{qq}\frac{\partial^2 W}{\partial x^2}+2D_{pq}\frac{\partial^2 W}{\partial \xi\partial x}.
	\end{aligned}
\end{equation}
Consequently, the 1D WPFP equation is given by:
\begin{align}
&\partial_t W
=-\xi\partial_{x}W-\Theta[V]W+D_{pp}\partial_{\xi\xi}W+D_{qq}\partial_{xx}W+2D_{pq}\partial_{x\xi}W+2\gamma\partial_{\xi}\left(\xi W\right), \label{eq:wpfp_1d} \\
&W(x,\xi,0)=W_0(x,\xi),\quad x,\xi\in\mathbb{R}. \nonumber 
\end{align}

Given the exponential decay of the solution in phase space, we truncate the infinite domain to a finite rectangle $[a,b]\times[c,d]$ and impose periodic boundary conditions. This approach reduces the original problem to the following initial-boundary value problem:
\begin{align}
&\partial_t W
=-\xi\partial_{x}W-\Theta[V]W+D_{pp}\partial_{\xi\xi}W+D_{qq}\partial_{xx}W+2D_{pq}\partial_{x\xi}W
+2\gamma\partial_{\xi}\left(\xi W\right),\nonumber\\
&\quad \quad\quad\quad \quad \quad\quad  x\in[a,b],\xi\in[c,d],  \label{eq:initial_boundary_problem} \\
&W(x,\xi,0)=W_0(x,\xi),\quad x\in[a,b],\xi\in[c,d], \nonumber \\
&W(a,\xi,t)=W(b,\xi,t),\quad\xi\in[c,d], \nonumber\\ 
&W(x,c,t)=W(x,d,t),\quad x\in[b,a], \nonumber
\end{align}
where $a,b,c,d \in \mathbb{R}$ are domain boundaries and $t \in (0,T],~T>0$ is the total simulation time.

For spatial discretization, we make a uniform partition of the computational domain $[a,b]\times[c,d]$ and define the following grid points:
\begin{align*}
x_j&=a+jh_x, \quad j=0,1,\cdots,M, \quad h_x=\frac{b-a}{M},\\ 
\xi_k&=c+kh_{\xi}, \quad k=0,1,\cdots,N, \quad h_{\xi}=\frac{d-c}{N}. 
\end{align*}
Here, $M$ and $N$ are even positive integers representing the number of grid points in $x$ and $\xi$ directions, respectively. For temporal discretization, we define time steps as 
$$
t_n=n\Delta t,\quad n=0,1,2,\cdots,P,\quad \Delta t=\frac{T}{P},
$$
where $P$ is a positive integer. We denote $W_{jk}^n\approx W(x_j,\xi_k,t_n)$ as the numerical approximation of the solution $W(x,\xi,t)$ at the grid point $(x_j,\xi_k,t_n)$.

\subsection{Time-splitting scheme}
\label{sec:time_splitting}

To numerically solve the WPFP equation \eqref{eq:initial_boundary_problem}, we employ a time-splitting strategy. This approach decomposes the original complex equation into several simpler equations, and these decomposed equations can be solved efficiently at each time step, respectively. 

Specifically, we rewrite Eq. \eqref{eq:initial_boundary_problem} as $\partial_t W=\mathcal{L}(W)$, where $\mathcal{L}(W)$ represents the right-hand side of \eqref{eq:initial_boundary_problem}.  The operator $\mathcal{L}$ can be decomposed into four parts: $\mathcal{L} = \mathcal{L}_1 + \mathcal{L}_2 + \mathcal{L}_3 + \mathcal{L}_4$, where each sub-operator $\mathcal{L}_i$ corresponds to a simpler subproblem:
\begin{align}
	\partial_t W=\mathcal{L}_1(W) &:= -\xi\partial_{x}W, \label{eq:split_L1} \\
	\partial_t W=\mathcal{L}_2(W) &:= -\Theta[V]W, \label{eq:split_L2} \\
	\partial_t W=\mathcal{L}_3(W) &:= D_{qq}\partial_{xx}W+2D_{pq}\partial_{x\xi}W+D_{pp}\partial_{\xi\xi}W, \label{eq:split_L3} \\
	\partial_t W=\mathcal{L}_4(W) &:= 2\gamma\partial_{\xi}\left(\xi W\right). \label{eq:split_L4}
\end{align}
The formal solution for $\partial_t W=\mathcal{L}(W)$ over a single time step $\Delta t$ can be expressed as $W(t_{n+1})=\exp(\Delta t\mathcal{L})W(t_n)$. To this end, we construct a numerical integrator based on the second-order Strang splitting scheme \cite{mclachlan2002splitting}, which is given by:
\begin{align*}
	\exp\left(\Delta t(A+B)\right)= \exp\left(\frac{\Delta t}{2} A\right)\exp\left(\Delta t B\right)\exp\left(\frac{\Delta t}{2} A\right)+O\left((\Delta t)^3\right).
\end{align*}

We apply this scheme in a nested fashion. First, we split the operator into $\mathcal{L} =\mathcal{L}_1+(\mathcal{L}_2 + \mathcal{L}_3+\mathcal{L}_4)$ and apply the Strang splitting scheme to it. Then, a further splitting of the second component leads to the final seven-stage scheme for advancing the solution from $t_n$ to $t_{n+1}$:
\begin{equation}
	\begin{aligned}
		W(t_{n+1}) \approx& \exp\left(\frac{\Delta t}{2}\mathcal{L}_1\right)\exp\left(\frac{\Delta t}{2}\mathcal{L}_2\right)
		\exp\left(\frac{\Delta t}{2}\mathcal{L}_3\right)\exp\left(\Delta t\mathcal{L}_4\right) \\
		&\exp\left(\frac{\Delta t}{2}\mathcal{L}_3\right)\exp\left(\frac{\Delta t}{2}\mathcal{L}_2\right)\exp\left(\frac{\Delta t}{2}\mathcal{L}_1\right)W(t_n).
	\end{aligned}
\end{equation}

Algorithmically, this means that for each time step, we solve the subproblems defined by \cref{eq:split_L1,eq:split_L2,eq:split_L3,eq:split_L4} in sequence, with the time increments prescribed by the exponents in the expression above.

For each subproblem $\partial_t W = \mathcal{L}_i(W)$, spatial discretization is achieved using the Fourier spectral method, chosen for its high accuracy and efficiency on periodic domains. The specific method applied to each problem is as follows:
\begin{itemize}
    \item The convection term~\eqref{eq:split_L1} is discretized using the Fourier spectral Galerkin method in the $x$-direction.
    \item The nonlocal term~\eqref{eq:split_L2} is handled with the Fourier spectral Galerkin method in the $\xi$-direction, coupled with a corresponding Fourier solver for the Poisson equation in the case of the Wigner-Poisson-Fokker-Planck system.
    \item The diffusion terms~\eqref{eq:split_L3} are discretized using the Fourier spectral Galerkin method in both the $x$- and $\xi$-directions.
    \item The friction term~\eqref{eq:split_L4} is primarily addressed using the Fourier spectral collocation method in the $\xi$-direction. A Fourier spectral Galerkin approach is also presented for comparison.
\end{itemize}
The detailed derivations for these discretizations are presented in the subsequent subsections.

\subsection{Discretization of the convection term~\eqref{eq:split_L1}}
\label{sec:convection_discretization}
This subsection details the numerical solution for the convection subproblem $\partial_t W = \mathcal{L}_1(W)$, as defined in \eqref{eq:split_L1}. This equation governing the transport of the solution in position space is discretized  by the Fourier spectral Galerkin method in the $x$-direction.

First, we approximate the solution $W(x,\xi,t)$ with a truncated Fourier series:
\begin{align}\label{eq:fourier_approx_x}
W(x,\xi,t)\approx\sum_{j=-M/2}^{M/2-1}\widehat{W}_j(\xi,t)\mathrm{e}^{\mathrm{i}\mu_j (x-a)},\quad\mu_j=\frac{2\pi j}{b-a},
\end{align}
where $\widehat{W}_j(\xi,t)$
are the semi-discrete Fourier coefficients obtained via the discrete Fourier transform (DFT):
\begin{align}\label{eq:DFT_IDFT}
\widehat{W}_j(\xi,t)\approx\frac{1}{b-a}\int_{a}^{b} W(x,\xi,t)\mathrm{e}^{-\mathrm{i}\mu_j(x-a)}\mathrm{d}x \approx \frac{1}{M}\sum_{m=0}^{M-1}W(x_m,\xi,t)\mathrm{e}^{-\mathrm{i}\mu_j(x_m-a)}, 
\end{align}
for $j=-M/2,\cdots,M/2-1.$

The Fourier-Galerkin method requires that the residual of the equation be orthogonal to each basis function. This leads to the Galerkin condition for each mode $j=-M/2,\cdots,M/2-1$:
\begin{align}\label{eq:galerkin_cond_x}
\left\langle\partial_t W-(-\xi\partial_x W),\mathrm{e}^{\mathrm{i}\mu_j(x-a)}\right\rangle_{x}=0,	 	
\end{align}
where $\langle f,g\rangle_{x}=\int_a^b f\overline{g}\mathrm{d}x, \overline{g}$ is the complex conjugate of $g$.
Substituting the Fourier expansion \eqref{eq:fourier_approx_x} into the Galerkin condition \eqref{eq:galerkin_cond_x} and using the orthogonality of the basis functions, we obtain a simple ordinary differential equation for each Fourier coefficient:
\begin{align*}
\partial_t\widehat{W}_j(\xi,t)
=-\mathrm{i}\mu_j\xi\widehat{W}_j(\xi,t).
\end{align*}
This ODE can be solved exactly, yielding the iteration scheme for the coefficients over a time step $\Delta t$:
\begin{align*}
\widehat{W}_j(\xi,t_{n+1})=\mathrm{e}^{-\mathrm{i}\mu_j\xi\Delta t}\widehat{W}_j(\xi,t_n).	
\end{align*}
Then $W(x,\xi,t)|_{t_{n+1}}$ is obtained by the inverse DFT (IDFT) of $\widehat{W}_j(\xi,t_{n+1})$. The implementation of this step is summarized in Algorithm~\ref{alg:step1_convection}.

\begin{algorithm}[H] 
\caption{Fourier spectral Galerkin discretization for solving \eqref{eq:split_L1}}\label{alg:step1_convection}
\begin{algorithmic}[1]
\State \textbf{Input} $W(x_m, \xi_k,t_n)$, for each $0\leq m \leq M-1,~0\leq k \leq N-1$.
\State \textbf{Step 1}  Compute the Fourier coefficients at $t_n$ by the 1D DFT in $x$:
    $$\widehat{W}_j(\xi_k,t_n) 
    =\frac{1}{M} \sum\limits_{m=0}^{M-1} W(x_m,\xi_k,t_n) \mathrm{e}^{-\mathrm{i}\mu_j(x_m-a)}.$$
\State \textbf{Step 2} Compute the Fourier coefficients at $t_{n+1}$: $$\widehat{W}_j(\xi_k,t_{n+1}) = \mathrm{e}^{-\mathrm{i}\mu_j\xi_k\Delta t}\widehat{W}_j(\xi_k,t_n).$$
\State \textbf{Step 3} Compute the solution at $t_{n+1}$ by the 1D IDFT:
    $$W(x_m,\xi_k,t_{n+1}) = \sum\limits_{j=-M/2}^{M/2-1}\widehat{W}_j(\xi_k,t_{n+1})\mathrm{e}^{\mathrm{i}\mu_j(x_m-a)}.$$
\State \textbf{Return} $W(x_m,\xi_k,t_{n+1})$ for each $0\leq m \leq M-1,~ 0\leq k \leq N-1$. 
\end{algorithmic}
\end{algorithm}

\subsection{Discretization of the nonlocal term \cref{eq:split_L2}}
\label{sec:nonlocal_discretization}

This subsection details the numerical scheme for the nonlocal subproblem $\partial_t W = \mathcal{L}_2(W)$. The discretization is based on a Fourier spectral Galerkin method in the momentum $\xi$- direction.

We approximate the solution $W(x,\xi,t)$ using the truncated Fourier series in $\xi$:
\begin{align}\label{eq:fourier_approx_xi}
W(x,\xi,t)\approx\sum_{k=-N/2}^{N/2-1}\widehat{W}_k(x,t)\mathrm{e}^{\mathrm{i}\nu_k (\xi-c)},\quad\nu_k=\frac{2\pi k}{d-c},
\end{align}
where the semi-discrete Fourier coefficients $\widehat{W}_k(x,t)$ are calculated via the Fourier transform
$$
\widehat{W}_k(x,t)=\frac{1}{d-c}\int_c^d W(x,\xi,t)\mathrm{e}^{-\mathrm{i}\nu_k(\xi-c)}\mathrm{d}\xi \approx \frac{1}{N}\sum_{l=0}^{N-1}W(x,\xi_l,t)\mathrm{e}^{-\mathrm{i}\nu_k(\xi_l-c)}.
$$
The Galerkin method requires the residual of the equation to be orthogonal to each basis function, leading to the condition:
\begin{align}\label{eq:galerkin_cond_xi_nonlocal}
\left\langle\partial_t W+\Theta[V]W,\mathrm{e}^{\mathrm{i}\nu_k(\xi-c)}\right\rangle_{\xi}=0,	 	
\end{align}
for each mode $k=-N/2,\cdots,N/2-1$, where $\langle f,g\rangle_{\xi}=\int_c^d f\overline{g}\,\mathrm{d}\xi$.

A key advantage of the Fourier basis is that the nonlocal operator $\Theta[V]$ becomes a simple multiplication in Fourier space. As detailed in \ref{appendix: nonloacl}, its action on the Fourier expansion is given by:
\begin{align}\label{eq:fourier_nonlocal_term}
-\Theta[V]W(x,\xi,t) = \sum_{k=-N/2}^{N/2-1} \delta V\left(x,\frac{\varepsilon\nu_k}{2},t\right)\widehat{W}_k(x,t)\mathrm{e}^{\mathrm{i}\nu_k(\xi-c)},
\end{align}
Substituting \eqref{eq:fourier_nonlocal_term} into the Galerkin condition~\eqref{eq:galerkin_cond_xi_nonlocal} yields a simple ordinary differential equation for each Fourier coefficient:
\begin{align*}
\partial_t\widehat{W}_k(x,t) = \delta V\left(x,\frac{\varepsilon\nu_k}{2},t\right)\widehat{W}_k(x,t).
\end{align*}
This ODE can be formally solved over a time step $\Delta t$. Using a first-order time-integration scheme, we obtain the iteration scheme:
\begin{align}\label{eq:update_fourier_nonlocal_term}
\widehat{W}_k(x,t_{n+1}) \approx \mathrm{e}^{\delta V(x,\frac{\varepsilon\nu_k}{2},t_n)\Delta t}\widehat{W}_k(x,t_n).	
\end{align}

The evaluation of the term $\delta V(x,\frac{\varepsilon\nu_k}{2},t_n)$ depends on the nature of the potential $V(x, t_n)$ distinguished the WFP equation from the WPFP equation.

For the WFP equation with a prescribed external potential, the procedure is straightforward. The potential $V(x,t_n)$ is known at each time step, allowing the term $\delta V$ to be computed directly from its definition.

For the WPFP equation with a self-consistent potential, the procedure is more computationally intensive and involves three main steps at each time point $t_n$:

\begin{enumerate}
    \item[\textbf{1.}] \textbf{Compute particle density.} The particle density $\rho(x,t_n)$ is obtained by integrating the Wigner function over the momentum coordinate:
    \[
    \rho(x,t_n) = \int_c^d W(x,\xi,t_n) \mathrm{d}\xi.
    \]

    \item[\textbf{2.}] \textbf{Solve the Poisson equation.}  Solving the potential $V(x,t_n)$ from the 1D Poisson equation:
    $$\partial_{xx}V(x,t_n)=\alpha \rho(x,t_n).$$ 
    To do this, we employ a Fourier spectral Galerkin method in the spatial direction. Firstly, approximating the potential and density by their truncated Fourier series:
    \begin{align}
        V(x,t)&\approx\sum_{j=-M/2}^{M/2-1}\widehat{V}_j(t)\mathrm{e}^{\mathrm{i}\mu_j(x-a)},\label{eq:fourier_approx_V}\\
        \rho(x,t)&\approx\sum_{j=-M/2}^{M/2-1}\widehat{\rho}_j(t)\mathrm{e}^{\mathrm{i}\mu_j(x-a)},\label{eq:fourier_approx_rho}
    \end{align}
    where $\mu_j=\frac{2\pi j}{b-a}$. Then enforcing the Galerkin condition
    \[
    \left\langle\partial_{xx}V(x,t)-\alpha\rho(x,t),\mathrm{e}^{\mathrm{i}\mu_j(x-a)}\right\rangle_{x}=0,
    \]
    on approximations \eqref{eq:fourier_approx_V} and \eqref{eq:fourier_approx_rho}, it transforms the Possion equation into a simple algebraic relation for the Fourier coefficients:
    \begin{align*}
    \widehat{V}_j(t_n)=-\frac{\alpha\widehat{\rho}_j(t_n)}{\mu_j^2} \quad \text{for } j\neq0, \quad \text{with } \widehat{V}_0(t_n)=0.
    \end{align*}

    \item[\textbf{3.}] \textbf{Construct the nonlocal term.} After we compute the potential's Fourier coefficients $\widehat{V}_j(t_n)$, the required term $\delta V$ for the scheme \eqref{eq:update_fourier_nonlocal_term} can be constructed using its Fourier series representation:
    \begin{align*}
    \delta V\left(x,\frac{\varepsilon\nu_k}{2},t_n\right) = \frac{\mathrm{i}}{\varepsilon}\sum_{j=-M/2}^{M/2-1}\widehat{V}_j(t_n)\mathrm{e}^{\mathrm{i}\mu_j(x-a)}\left[\mathrm{e}^{\mathrm{i}\frac{\varepsilon\mu_j\nu_k}{2}}-\mathrm{e}^{-\mathrm{i}\frac{\varepsilon\mu_j\nu_k}{2}}\right].
    \end{align*}
\end{enumerate}

The entire procedure for solving subproblem~\eqref{eq:split_L2} for both cases is summarized in Algorithm~\ref{alg:step2_nonlocal}.

\begin{algorithm}[H]
\caption{Fourier spectral Galerkin discretization for solving \eqref{eq:split_L2}}\label{alg:step2_nonlocal}
\begin{algorithmic}[1]
\State \textbf{Input} $W(x_m, \xi_l,t_n)$, for each $0\leq m \leq M-1,~0\leq l \leq N-1$.
    \State \textbf{Step 1}  Compute the potential $\rho(x_m,t_n)$ by
    $$\rho(x_m,t_n) = \int_c^d W(x_m,\xi,t_n) \mathrm{d}\xi \approx \sum_{l=0}^{N-1} W(x_m,\xi_l,t_n)h_{\xi},$$
    and the Fourier coefficients $\widehat{\rho}_j(t_n)$ by the DFT in $x$: 
    $$\widehat{\rho}_j(t_n)=\frac{1}{M} \sum\limits_{m=0}^{M-1} \rho(x_m,t_n) \mathrm{e}^{-\mathrm{i}\mu_j(x_m-a)}.$$ 
    \State \textbf{Step 2} Compute 
    the Fourier coefficients $\widehat{V}_j(t_n)$ by
    $$\widehat{V}_0(t_n)=0, \widehat{V}_j(t_n)=-\frac{\alpha\widehat{\rho}_j(t_n)}{\mu_j^2},j\neq0,$$
    and $\delta V(x_m,\frac{\varepsilon\nu_k}{2},t_n)$ by the IDFT:
    $$\delta V(x_m,\frac{\varepsilon\nu_k}{2},t_n)=\frac{\mathrm{i}}{\varepsilon}\sum_{j=-M/2}^{M/2-1}\widehat{V}_j(t_n)\mathrm{e}^{\mathrm{i}\mu_j(x_m-a)}\left[\mathrm{e}^{\mathrm{i}\frac{\varepsilon\mu_j\nu_k}{2}}-\mathrm{e}^{-\mathrm{i}\frac{\varepsilon\mu_j\nu_k}{2}}\right].$$
    \State \textbf{Step 3} (For the WFP case, we can skip Step 1 and 2 and go to Step 3 directly.) Compute the Fourier coefficients at $t_n$ by the 1D DFT in $\xi$:
    $$\widehat{W}_k(x_m,t_{n})= \frac{1}{N}\sum_{l=0}^{N-1}W(x_m,\xi_l,t_n)\mathrm{e}^{-\mathrm{i}\nu_k(\xi_l-c)},$$
    then the coefficients at $t_{n+1}$: 
             $$\widehat{W}_k(x_m,t_{n+1})=\mathrm{e}^{\delta V(x_m,\frac{\varepsilon\nu_k}{2},t_n)\Delta t}\widehat{W}_k(x_m,t_n).$$  
     \State \textbf{Step 4}  Compute the solution $W(x_m,\xi_l,t_{n+1})$ by the 1D IDFT of $\widehat{W}_k(x_m,t_{n+1}).$
\State \textbf{Return} $W(x_m,\xi_l,t_{n+1})$  for each $0\leq m \leq M-1,~ 0\leq l \leq N-1$. 
\end{algorithmic}
\end{algorithm}

\subsection{Discretization of the diffusion terms \eqref{eq:split_L3}}
\label{sec:diffusion_discretization}
This subsection presents the numerical scheme for the diffusion subproblem $\partial_t W = \mathcal{L}_3(W)$ defined by \eqref{eq:split_L3}. This linear partial differential equation, which encompasses diffusion in both position and momentum as well as cross-diffusion, is ideally suited to be solved by the Fourier spectral Galerkin method applied in both the $x$- and $\xi$- directions.

We approximate the solution $W(x,\xi,t)$ in phase space using a two-dimensional truncated Fourier series:
\begin{align}\label{eq:fourier_approx_phase_space}
	W(x,\xi,t)\approx\sum_{j=-M/2}^{M/2-1}\sum_{k=-N/2}^{N/2-1}\widehat{W}_{jk}(t)\mathrm{e}^{\mathrm{i}\mu_j(x-a)+\mathrm{i}\nu_k(\xi-c)},
\end{align}
where $\mu_j=\frac{2\pi j}{b-a}$, $\nu_k=\frac{2\pi k}{d-c}$ and the two-dimensional Fourier coefficients $\widehat{W}_{jk}(t)$ are calculated by
\begin{align*}
	\widehat{W}_{jk}(t)&\approx\frac{1}{(b-a)(d-c)}\int_a^b\int_c^d W(x,\xi,t)\mathrm{e}^{-\mathrm{i}\mu_j(x-a)-\mathrm{i}\nu_k(\xi-c)}\mathrm{d}\xi\mathrm{d}x \\
	&\approx\frac{1}{MN}\sum_{m=0}^{M-1}\sum_{l=0}^{N-1}W(x_m,\xi_l,t)\mathrm{e}^{-\mathrm{i}\mu_j(x_m-a)-\mathrm{i}\nu_k(\xi_l-c)}.	
\end{align*}
The Fourier-Galerkin method requires that the residual of the equation be orthogonal to each basis function. This yields the following Galerkin condition for each mode $(j,k)$:
\begin{align}\label{eq:galerkin_cond_phase_space}
	\left\langle\partial_t W-\left(D_{qq}\partial_{xx}W+2D_{pq}\partial_{x\xi}W+D_{pp}\partial_{\xi\xi}W\right),\mathrm{e}^{\mathrm{i}\mu_j(x-a)+\mathrm{i}\nu_k(\xi-c)}\right\rangle_{x,\xi}= 0,	
\end{align}
where the inner product is defined over the phase space as $\langle f,g\rangle_{x,\xi}=\int_a^b\int_c^d f\overline{g}\,\mathrm{d}x\mathrm{d}\xi$.

Substituting the Fourier expansion \eqref{eq:fourier_approx_phase_space} into the Galerkin condition \eqref{eq:galerkin_cond_phase_space} leads to a simple, uncoupled ordinary differential equation for each Fourier coefficient:
\begin{align*}
	\partial_{t}\widehat{W}_{jk}(t)=\left(-D_{qq}\mu_j^2-2D_{pq}\mu_j\nu_k-D_{pp}\nu_k^2\right)\widehat{W}_{jk}(t).
\end{align*}
This ODE is linear with constant coefficients and has an exact solution expresed by the following explicit update rule over a time step $\Delta t$:
\begin{align}
	\widehat{W}_{jk}(t_{n+1})=\mathrm{e}^{\left(-D_{qq}\mu_j^2-2D_{pq}\mu_j\nu_k-D_{pp}\nu_k^2\right)\Delta t}\widehat{W}_{jk}(t_{n}).	
\end{align}
The practical implementation of this method is summarized in Algorithm~\ref{alg:step3_diffusion}.

\begin{algorithm}[H] 
\caption{Fourier spectral Galerkin discretization for solving \eqref{eq:split_L3}}\label{alg:step3_diffusion}
\begin{algorithmic}[1]
\State \textbf{Input} $W(x_m, \xi_k,t_n)$, for each $0\leq m \leq M-1,~0\leq k \leq N-1$.
\State \textbf{Step 1}  Compute the Fourier coefficients at $t_n$ by the 2D DFT in $x$ and $\xi$:
    $$\widehat{W}_{jk}(t_n) 
    =\frac{1}{MN}\sum_{m=0}^{M-1}\sum_{l=0}^{N-1}W(x_m,\xi_l,t)\mathrm{e}^{-\mathrm{i}\mu_j(x_m-a)-\mathrm{i}\nu_k(\xi_l-c)}.$$
\State \textbf{Step 2} Compute the Fourier coefficients at $t_{n+1}$: $$\widehat{W}_{jk}(t_{n+1}) = \mathrm{e}^{\left(-D_{qq}\mu_j^2-2D_{pq}\mu_j\nu_k-D_{pp}\nu_k^2\right)\Delta t}\widehat{W}_{jk}(t_{n}).$$
\State \textbf{Step 3} Compute the solution at $t_{n+1}$ by the 2D IDFT:
    $$W(x_m,\xi_k,t_{n+1}) = \sum_{j=-M/2}^{M/2-1}\sum_{k=-N/2}^{N/2-1}\widehat{W}_{jk}(t_{n+1})\mathrm{e}^{\mathrm{i}\mu_j(x-a)+\mathrm{i}\nu_k(\xi-c)}.$$
\State \textbf{Return} $W(x_m,\xi_k,t_{n+1})$ for each $0\leq m \leq M-1,~ 0\leq k \leq N-1$. 
\end{algorithmic}
\end{algorithm}

\subsection{Discretization of the friction term \eqref{eq:split_L4}}
\label{sec:friction_collocation_discretization}
This subsection presents two numerical approaches for solving the friction subproblem $\partial_t W = \mathcal{L}_4(W)$ given by \eqref{eq:split_L4}. The primary scheme is based on the Fourier spectral collocation method, while the Fourier spectral Galerkin method is also implemented for comparison.

\subsubsection{The Fourier spectral collocation method}
We approximate the solution $W(x,\xi,t)$ using the spectral collocation expansion with respect to $\xi$:
\begin{align}\label{eq: spectral collocation approximation}
W(x,\xi,t)=\sum_{j=0}^{N-1}W(x,\xi_j,t)H_j(\xi),	
\end{align}
where $W(x,\xi_j,t)$ are the solution values at the collocation points $\xi_j$. The basis functions, $H_j(\xi)$, are constructed to satisfy $H_j(\xi_k)=\delta_{jk}$, ensuring that the approximation exactly matches the values $W(x,\xi_k,t)$ at the nodes. These basis functions can be constructed from the symmetric Fourier spectral approximation formula:
\begin{align*}
W(x,\xi,t)
&=\sum_{k=-\frac{N}{2}}^{\frac{N}{2}}\frac{1}{c_k}\widehat{W}_k(x,t)\mathrm{e}^{\mathrm{i}\mu_k (\xi-c)}\approx\sum_{k=-\frac{N}{2}}^{\frac{N}{2}}\frac{1}{c_k}\left[\frac{1}{N}\sum_{j=0}^{N-1}W(x,\xi_j,t)\mathrm{e}^{-\mathrm{i}\mu_k(\xi_j-c)}\right]\mathrm{e}^{\mathrm{i}\mu_k (\xi-c)}\\
&=\sum_{j=0}^{N-1}W(x,\xi_j,t)\left[\frac{1}{N}\sum_{k=-\frac{N}{2}}^{\frac{N}{2}}\frac{1}{c_k}\mathrm{e}^{\mathrm{i}\mu_k(\xi-\xi_j)}\right]\\
&=\sum_{j=0}^{N-1}W(x,\xi_j,t)\frac{1}{N}\sin\left(\frac{N\pi}{d-c}(\xi-\xi_j)\right)\cot\left(\frac{\pi}{d-c}(\xi-\xi_j)\right).	
\end{align*}
where $c_k$ is the modified symmetry Fourier spectral approximation coefficient which satisfies
\begin{align*}
c_k=\begin{cases}
1, |k|<\frac{N}{2},\\
2, |k|=\frac{N}{2}.	
\end{cases}	
\end{align*}
Therefore, the explicit form of the basis function $H_j(\xi)$ is
\begin{align}\label{eq: }
H_j(\xi)=\frac{1}{N}\sin\left(\frac{N\pi}{d-c}(\xi-\xi_j)\right)\cot\left(\frac{\pi}{d-c}(\xi-\xi_j)\right).	
\end{align}

The core principle of the collocation method is to enforce that the equation holds exactly at each grid point $\xi_k$, for $k=0,1,\cdots,N-1$. Thus, the collocation condition is given as:
\begin{align}\label{eq: collocation condition}
	\partial_t W\big|_{\xi=\xi_k}=2\gamma\partial_{\xi}(\xi W)\big|_{\xi=\xi_k},\quad k=0,1,2,\cdots,N-1.
\end{align}
Substituting \cref{eq: spectral collocation approximation} into condition \cref{eq: collocation condition} and using the property $H_j(\xi_k)=\delta_{jk}$, we have
\begin{align}\label{eq: collocation method-equation}
\partial_t W(x,\xi_k,t)=2\gamma W(x,\xi_k,t)+2\gamma\xi_k\sum_{j=0}^{N-1}d_{kj} W(x,\xi_j,t),	
\end{align}
where the entries $d_{kj}:=H_j^{\prime}(\xi_k)$ form the spectral differentiation matrix, whose formula is
\begin{align}\label{eq: derivative of H_j}
d_{kj}
=\begin{cases}
(-1)^{k+j}\frac{\pi}{d-c}\cot\left(\frac{\pi}{N}(k-j)\right),&k\neq j\\
0,&k=j.	
\end{cases}
\end{align}

To express this system in a more compact form, we define the solution matrix $\bm{W}$, the diagonal coordinate matrix $\bm{\Lambda}$, and the spectral differentiation matrix $\bm{D}$ as:
\begin{align*}
\bm{W}(t)&=
\begin{bmatrix}
W(x_0,\xi_0,t) & W(x_1,\xi_0,t) & \cdots & W(x_{N-1},\xi_0,t) \\
W(x_0,\xi_1,t) & W(x_1,\xi_1,t) & \cdots & W(x_{N-1},\xi_{1},t) \\
\vdots & \vdots & \ddots & \vdots \\
W(x_0,\xi_{N-1},t) & W(x_1,\xi_{N-1},t) & \cdots & W(x_{N-1},\xi_{N-1},t)
\end{bmatrix},\\
\bm{\Lambda}&=\begin{bmatrix}
\xi_0 & 0   & \cdots & 0 \\
0   & \xi_1 & \cdots & 0 \\
\vdots & \vdots & \ddots & \vdots \\
0   & 0   & \cdots & \xi_{N-1}
\end{bmatrix},
\bm{D}=\begin{bmatrix}
d_{00} & d_{01} & \cdots & d_{0,N-1} \\
d_{10} & d_{11} & \cdots & d_{1,N-1} \\
\vdots & \vdots & \ddots & \vdots \\
d_{N-1,0} & d_{N-1,1} & \cdots & d_{N-1,N-1}
\end{bmatrix}.	
\end{align*}
The ODE system \cref{eq: collocation method-equation} can now be written in matrix form as:
\begin{align*}
\partial_t\bm{W}=2\gamma(\bm{I}+\bm{\Lambda}\bm{D})\bm{W}.	
\end{align*}
where $\bm{I}$ is the identity matrix. Solving this ODE system from time $t_n$ to $t_{n+1}$ yields
\begin{align*}
\bm{W}(t_{n+1})=\exp\left[2\gamma\left(\bm{I}+\bm{\Lambda}\bm{D}\right)\Delta t\right]\bm{W}(t_n).	
\end{align*}
The complete procedure is summarized in \cref{alg:step4_collocation}.

\begin{algorithm}[H]
\caption{Fourier spectral collocation discretization for solving \eqref{eq:split_L4}}\label{alg:step4_collocation}
\begin{algorithmic}[1]
\State \textbf{Input} The solution matrix $\bm{W}(t_n)=[W(x_m, \xi_k,t_n)]_{m,k}$, for each $0\leq m \leq M-1,~0\leq k \leq N-1$.
\State \textbf{Step 1} Precompute the differential matrix $\bm{D}=[d_{kl}]_{k,l=0}^{N-1}$ by \eqref{eq: derivative of H_j}, the diagonal matrix $\mathbf{\Lambda} = \operatorname{diag}(\xi_0, \dots, \xi_{N-1})$, and then the iteration matrix 
$$\mathbf{A} = 2\gamma(\mathbf{I}+\mathbf{\Lambda}\mathbf{D}).$$
\State \textbf{Step 2} Compute the solution matrix at $t_{n+1}$:
$$\mathbf{W}(t_{n+1})= \exp(A\Delta t)\mathbf{W}(t_n).$$  
\State \textbf{Return} $\bm{W}(t_{n+1})$ for each $0\leq m \leq M-1,~0\leq k \leq N-1$.
\end{algorithmic}
\end{algorithm}

\subsubsection{The Fourier spectral Galerkin method}
\label{sec:friction_galerkin_discretization}
An alternative approach to solving Eq.~\eqref{eq:split_L4} is the Fourier spectral Galerkin method. The method begins by approximating the solution $W$ with its Fourier series~\eqref{eq:fourier_approx_xi} and applying the Galerkin condition, which requires the residual to be orthogonal to the basis functions:
\begin{align}\label{eq:galerkin_cond_xi_friction}
\left\langle\partial_t W-\left(2\gamma\partial_{\xi}(\xi W)\right),\mathrm{e}^{\mathrm{i}\nu_k(\xi-c)}\right\rangle_{\xi}=0,	 	
\end{align}
for each mode $k=-N/2,\cdots,N/2-1$. Then substituting the Fourier approximation into~\eqref{eq:galerkin_cond_xi_friction} and using integration by parts contribute to a system of coupled ODEs for the Fourier coefficients:
\begin{align}\label{eq:galerkin_method_equation_friction_final}
\partial_t\widehat{W}_k(x,t)=\frac{2\gamma}{d-c}\left[(d-c+\mathrm{i}\pi k(d+c))\widehat{W}_k(x,t)+(d-c)\sum_{l=-N/2,l\neq k}^{N/2-1}\frac{l}{l-k}\widehat{W}_l(x,t)\right].
\end{align}

This system can be expressed in a concise matrix form. Let us define the matrix of Fourier coefficients $\widehat{\bm{W}}(t)$ and the operational matrices $\bm{E}$ and $\bm{F}$ as follows:
\begin{align*}
\widehat{\bm{W}}(t)&=
\begin{bmatrix}
\widehat{W}_0(x_0,t) & \widehat{W}_0(x_1,t) & \cdots & \widehat{W}_0(x_{N-1},t) \\
\widehat{W}_1(x_0,t) & \widehat{W}_1(x_1,t) & \cdots & \widehat{W}_1(x_{N-1},t) \\
\vdots & \vdots & \ddots & \vdots \\
\widehat{W}_{\frac{N}{2}-1}(x_0,t) & \widehat{W}_{\frac{N}{2}-1}(x_1,t) & \cdots & \widehat{W}_{\frac{N}{2}-1}(x_{N-1},t) \\
\widehat{W}_{-\frac{N}{2}}(x_0,t) & \widehat{W}_{-\frac{N}{2}}(x_1,t) & \cdots & \widehat{W}_{-\frac{N}{2}}(x_{N-1},t) \\
\vdots & \vdots & \ddots & \vdots \\
\widehat{W}_{-1}(x_0,t) & \widehat{W}_{-1}(x_1,t) & \cdots & \widehat{W}_{-1}(x_{N-1},t)
\end{bmatrix},\\
\bm{E}&=\mathrm{i}\pi\frac{d+c}{d-c}\operatorname{diag}\left(0,1,\cdots,\frac{N}{2}-1,-\frac{N}{2},\cdots,-1\right),
\bm{F}=(f_{kl}),
\end{align*}
where
\begin{align*}
f_{kl}=\cdot\begin{cases}
\frac{l}{l-k},\quad &l\neq k,\\
0, \quad &l=k.	
\end{cases}	
\end{align*}
Using these definitions, we can write \cref{eq:galerkin_method_equation_friction_final} in matrix form:
\begin{align*}
\partial_t\widehat{\bm{W}}(t)=2\gamma(\bm{I}+\bm{E}+\bm{F})\widehat{\bm{W}}(t).	
\end{align*}
where $\bm{I}$ is the identity matrix. Integrating this ODE from $t_n$ to $t_{n+1}$ gives the solution at the next time step:
\begin{align*}
\widehat{\bm{W}}(t_{n+1})=\exp\left[2\gamma\left(\bm{I}+\bm{E}+\bm{F}\right)\Delta t\right]\widehat{\bm{W}}(t_n).	
\end{align*}

\subsection{Overall Time-Splitting Fourier pseudospectral algorithm}
\label{sec:overall_algorithm}
The numerical methods for solving the sub-problems \cref{eq:split_L1,eq:split_L2,eq:split_L3,eq:split_L4} were detailed in the preceding subsections. By combining these steps, we arrive at the complete numerical algorithm for the 1D WPFP system, which is summarized in Algorithm \ref{alg: total}.

\begin{algorithm}[H]
\caption{Time-Splitting Fourier pseudospectral algorithm for the 1D WPFP Equation}
\label{alg: total}
\begin{algorithmic}[1]
\State \textbf{Input:} 
\Statex \quad Initial condition $W_0(x, \xi)$.
\Statex \quad Computational domain $[a,b]\times[c,d]$ with an $M \times N$ grid.
\Statex \quad Final simulation time $T$ and time step $\Delta t$.
\Statex \quad All relevant physical and model parameters.
\State \textbf{Output:} The final numerical solution $W(x, \xi, T)$.

\vspace{0.5em}
\State \textbf{Initialization:}
\State Evaluate the initial condition on the grid: $W^0 \gets W_0(x_j, \xi_k)$.
\State Compute the total number of time steps: $P \gets \lfloor T / \Delta t \rfloor$.

\vspace{0.5em}
\For{$n = 0$ to $P - 1$}
    \State Let $W_{\text{temp}} \gets W^n$.
    
    \State Apply the convection step (Alg.~\ref{alg:step1_convection}) to $W_{\text{temp}}$ with time increment $\Delta t/2$ and update $W_{\text{temp}}$.
    \State Apply the nonlocal step (Alg.~\ref{alg:step2_nonlocal}) to $W_{\text{temp}}$ with time increment $\Delta t/2$ and update $W_{\text{temp}}$.
    \State Apply the diffusion step (Alg.~\ref{alg:step3_diffusion}) to $W_{\text{temp}}$ with time increment $\Delta t/2$ and update $W_{\text{temp}}$.
    \State Apply the friction step (Alg.~\ref{alg:step4_collocation}) to $W_{\text{temp}}$ with time increment $\Delta t$ and update $W_{\text{temp}}$.
    \State Apply the diffusion step (Alg.~\ref{alg:step3_diffusion}) to $W_{\text{temp}}$ with time increment $\Delta t/2$ and update $W_{\text{temp}}$.
    \State Apply the nonlocal step (Alg.~\ref{alg:step2_nonlocal}) to $W_{\text{temp}}$ with time increment $\Delta t/2$ and update $W_{\text{temp}}$.
    \State Apply the convection step (Alg.~\ref{alg:step1_convection}) to $W_{\text{temp}}$ with time increment $\Delta t/2$ and update $W_{\text{temp}}$.
    
    \State Set $W^{n+1} \gets W_{\text{temp}}$.
\EndFor

\vspace{0.5em}
\State \Return $W^{P}$.
\end{algorithmic}
\end{algorithm}

\section{Numerical Results}
\label{sec: numerical results}

In this section, we present numerical results to validate the proposed time-splitting Fourier spectral (TSSP) method. First, we verify the proposed method's spectral accuracy in both the spatial $x$- and momentum $\xi$- directions, as well as its second-order accuracy in time. These convergence tests are performed for the WFP system with both harmonic and non-harmonic potentials, and for the WPFP system, respectively. Second, we numerically investigate the existence of steady states for the WFP system with harmonic potentials containing a small perturbation or a large perturbation.
Finally, we examine the convergence to a steady state for the WPFP system.

In all subsequent examples, the initial Wigner function is defined as a Gaussian wavepacket in phase space centered at $(x_0, \xi_0)$:
\begin{align*}
f(x,\xi,0)=\frac{\sqrt{a_{11}a_{22}-a_{12}^2}}{\pi\varepsilon}\exp{\left(-\frac{a_{11}(x-x_0)^2+a_{22}(\xi-\xi_0)^2+2a_{12}(x-x_0)(\xi-\xi_0)}{\varepsilon}\right)},
\end{align*}
where the coefficients $a_{11}, a_{12}, a_{22}$ form the covariance matrix of the wavepacket, determining its shape, size, and orientation in phase space.

\subsection{Convergence Test}

\textbf{Example 1: Harmonic Potential in WFP Equation}

We begin by considering a harmonic potential $V(x)=\frac{1}{2}x^2+x$. For this potential, an analytical solution based on wavepacket dynamics (solving by the method of characteristics) exists for the 1D WFP equation. This analytical solution serves as the reference against which we benchmark our numerical results. The model parameters are set to $\varepsilon=0.1$, $D_{pp}=0.2$, $D_{qq}=0.2$, $D_{pq}=0.05$, and $\gamma=1$. The initial conditions are defined by $a_{11}=a_{22}=-1$, $x_0=0.1$, and $\xi_0=-0.2$. The computational domain is $[-2,2]\times[-2,2]$, and the simulation is run up to a final time of $T=0.5$. The reference solution is computed on the same spatial mesh using a fine time step of $\Delta t=2^{-10}$.

\begin{figure}[htbp]
    \centering
    \begin{subfigure}{0.48\textwidth}
        \centering
        \includegraphics[width=\textwidth]{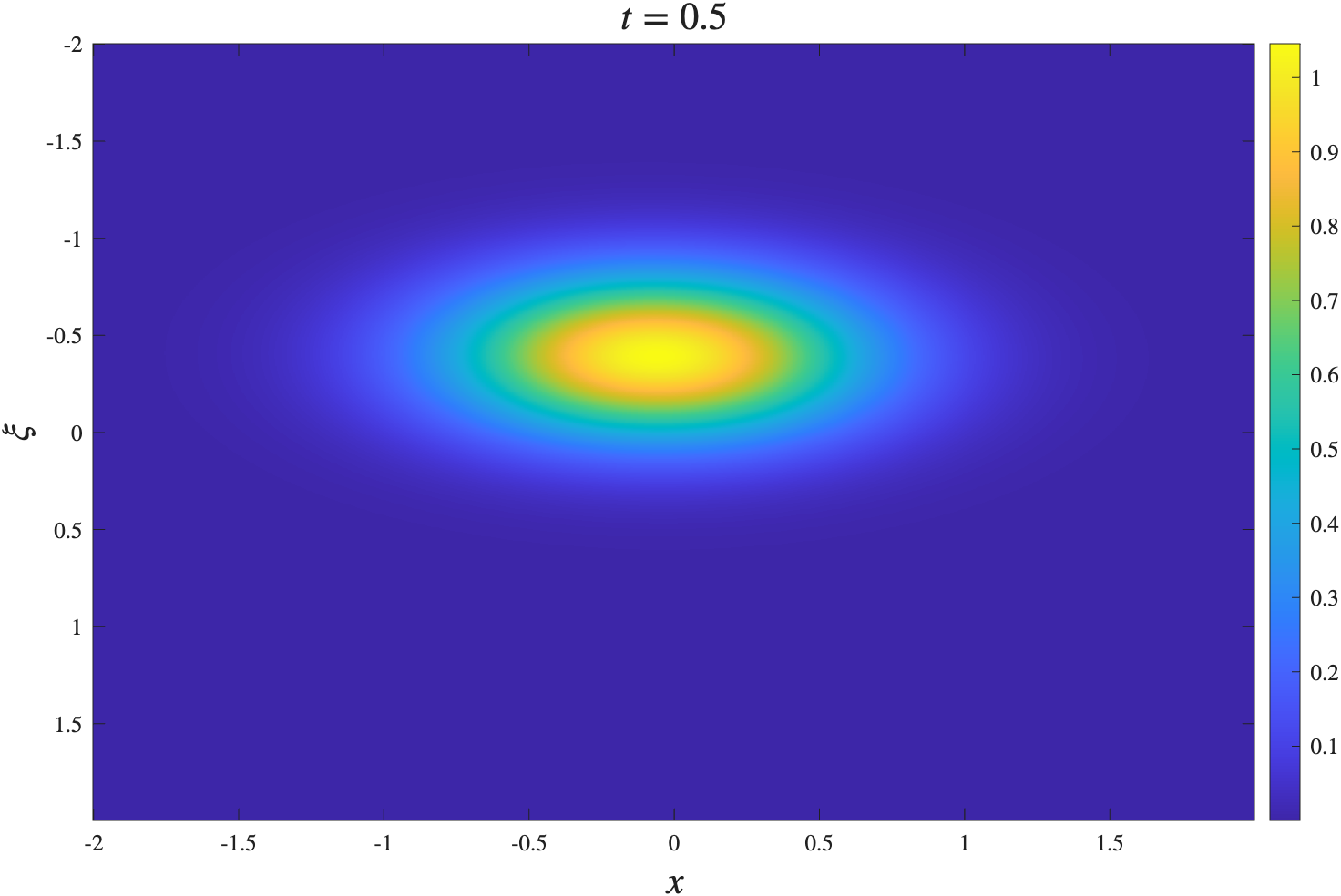}
        \label{fig: ep1-exact}
    \end{subfigure}
    \begin{subfigure}{0.48\textwidth}
        \centering
        \includegraphics[width=\textwidth]{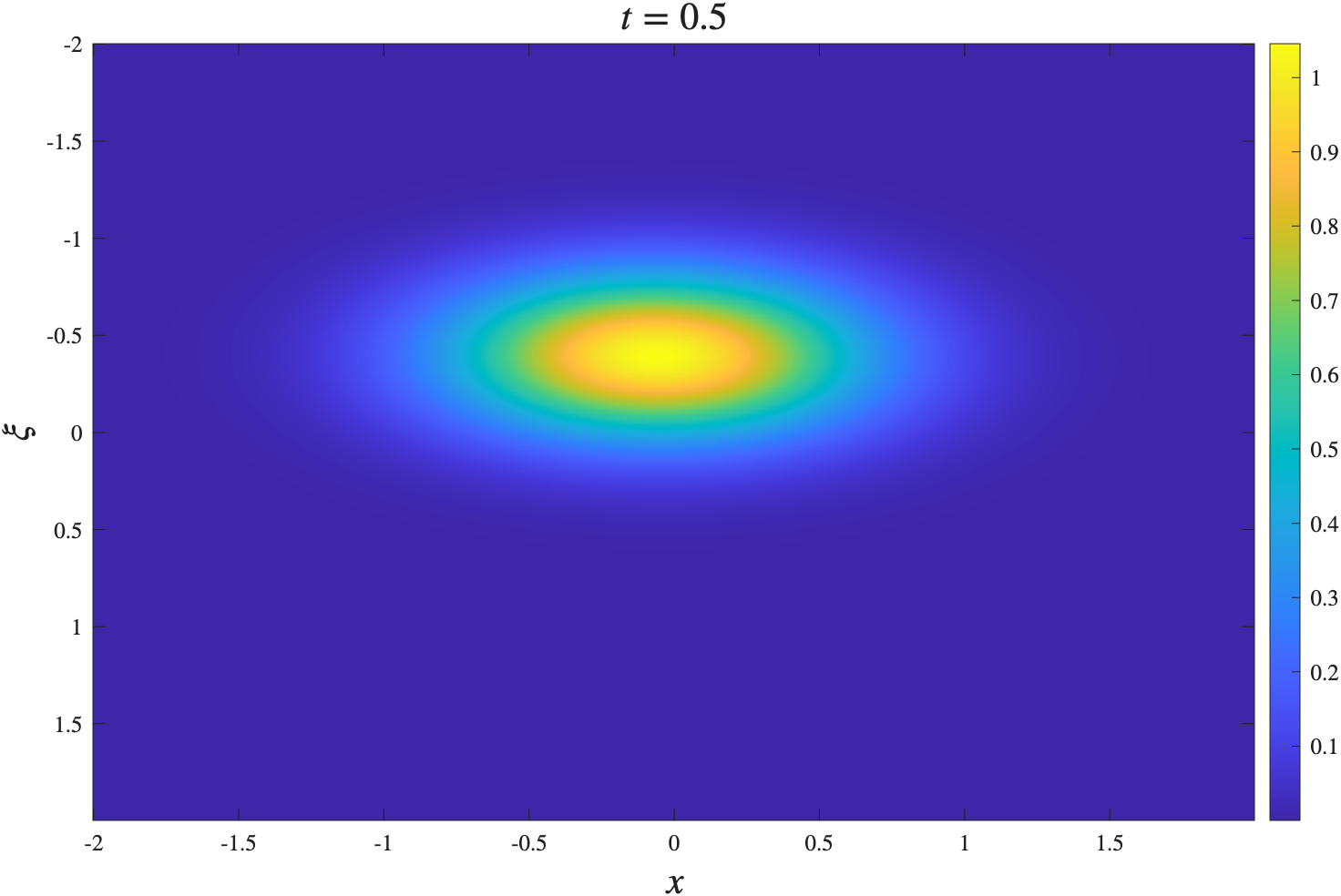}
        \label{fig: ep1-tssp}
    \end{subfigure}
    \caption{Heatmaps of the Wigner function at $T=0.5$ for Example 1 with semiclassical parameter $\varepsilon=0.1$. Left: exact solution by the method of characteristics; Right: numerical solution by the time-splitting spectral method.}
    \label{fig: example 1-harmonic potential}
\end{figure}

\begin{figure}[htbp]
    \centering
    \begin{subfigure}{0.32\textwidth}
        \centering
        \includegraphics[width=\textwidth]{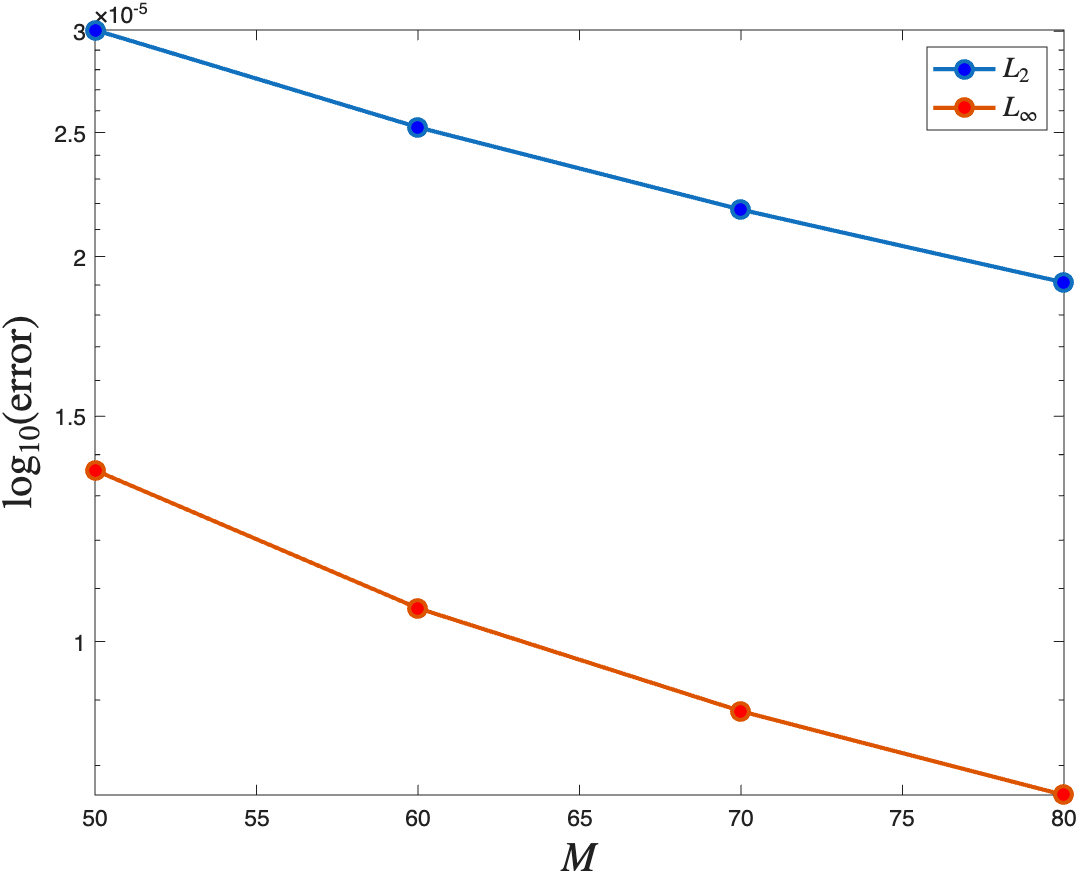}
        \label{fig: ep1-error-x}
    \end{subfigure}
    \begin{subfigure}{0.32\textwidth}
        \centering
        \includegraphics[width=\textwidth]{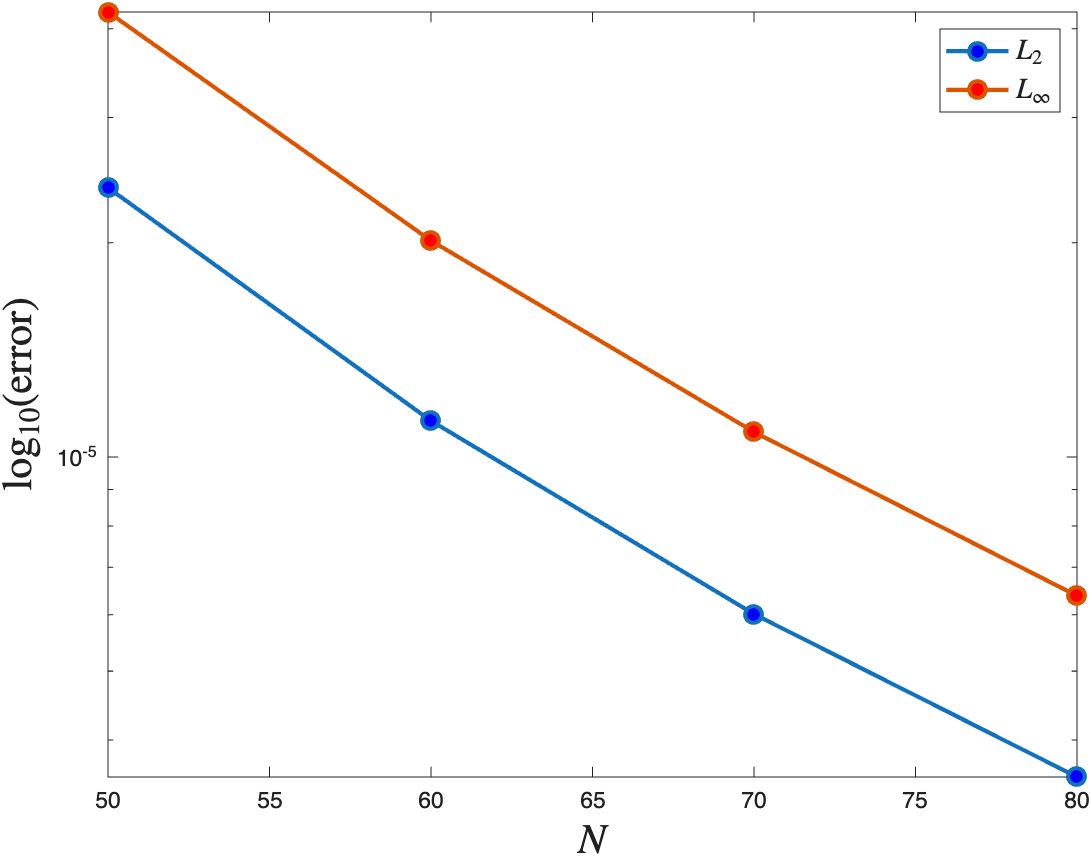}
        \label{fig: ep1-error-xi}
    \end{subfigure}
    \begin{subfigure}{0.32\textwidth}
        \centering
        \includegraphics[width=\textwidth]{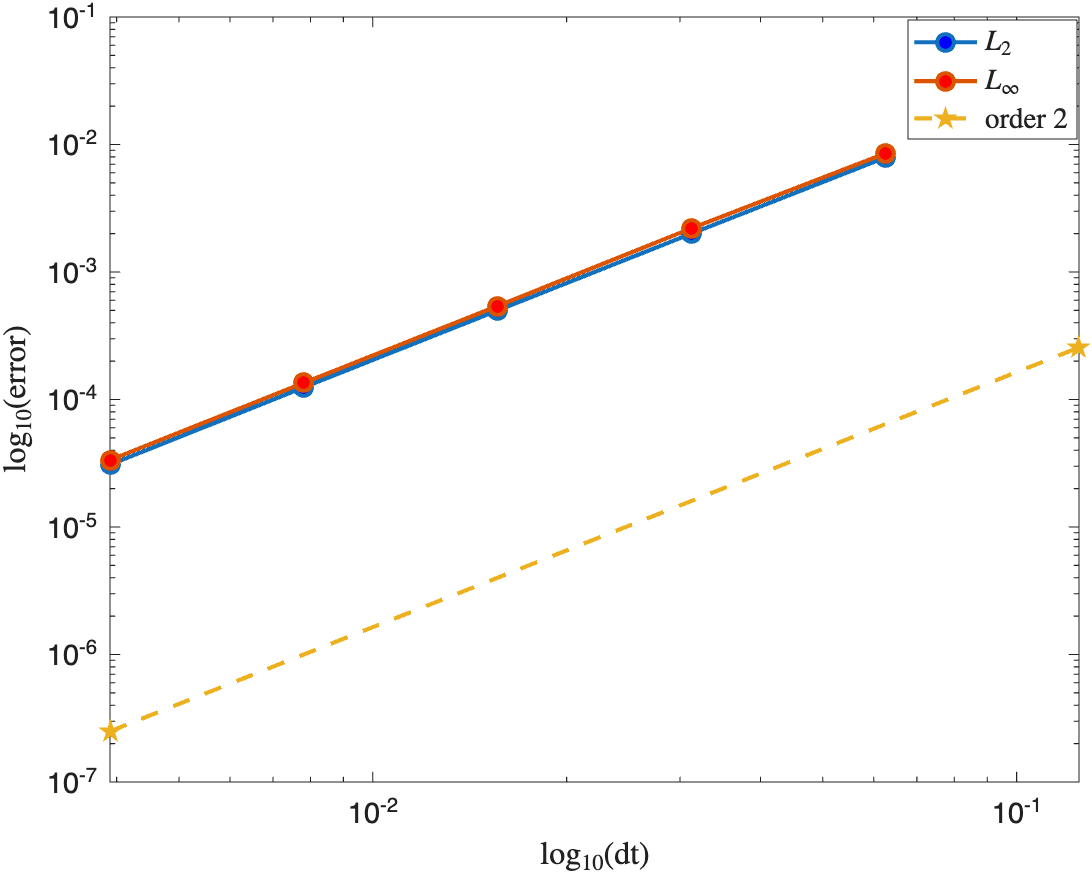}
        \label{fig: ep1-error-time}
    \end{subfigure}
    \caption{Convergence rates for \textbf{Example 1} at $T=0.5$ with respect to $M$ (left), $N$ (middle), and time step $\Delta t$ (right). The red and blue lines correspond to the errors in $L^2$ and $L^{\infty}$ norm, respectively.}
    \label{fig: example 1-harmonic potential error order}
\end{figure}

Figure \ref{fig: example 1-harmonic potential} comparing the heatmaps of the Wigner function computed by the analytical solution and our numerical method, demonstrates that the proposed method accurately captures the system's dynamical behavior. To evaluate the convergence rates, we plot the errors performed in different refined meshegrids in Figure \ref{fig: example 1-harmonic potential error order}.
From this figure, we confirm the spectral convergence of the TSSP method with respect to the number of grid points $M$ and $N$. The right panel of the figure plots the error against the time step $\Delta t$ (with fixed spatial resolution $N=M=2^9$), confirming the expected second-order temporal convergence. Similar results were also obtained through the TSSP-Galerkin method (solving the friction term by Galerkin method) which are omitted here for brevity.

\textbf{Example 2: Non-harmonic Potential in WFP Equation}

We now consider a non-harmonic potential, $V(x)=(x^2-1)^2$, to test the method's performance in a more complex scenario. All other parameters and initial conditions remain identical to those in Example 1. Since an analytical solution is unavailable for this potential, we use a highly resolved numerical simulation as the reference solution, computed with a spatial grid of $2^{10} \times 2^{10}$ points and a time step of $\Delta t=2^{-10}$.

\begin{figure}[htbp]
    \centering
    \begin{subfigure}{0.32\textwidth}
        \centering
        \includegraphics[width=\textwidth]{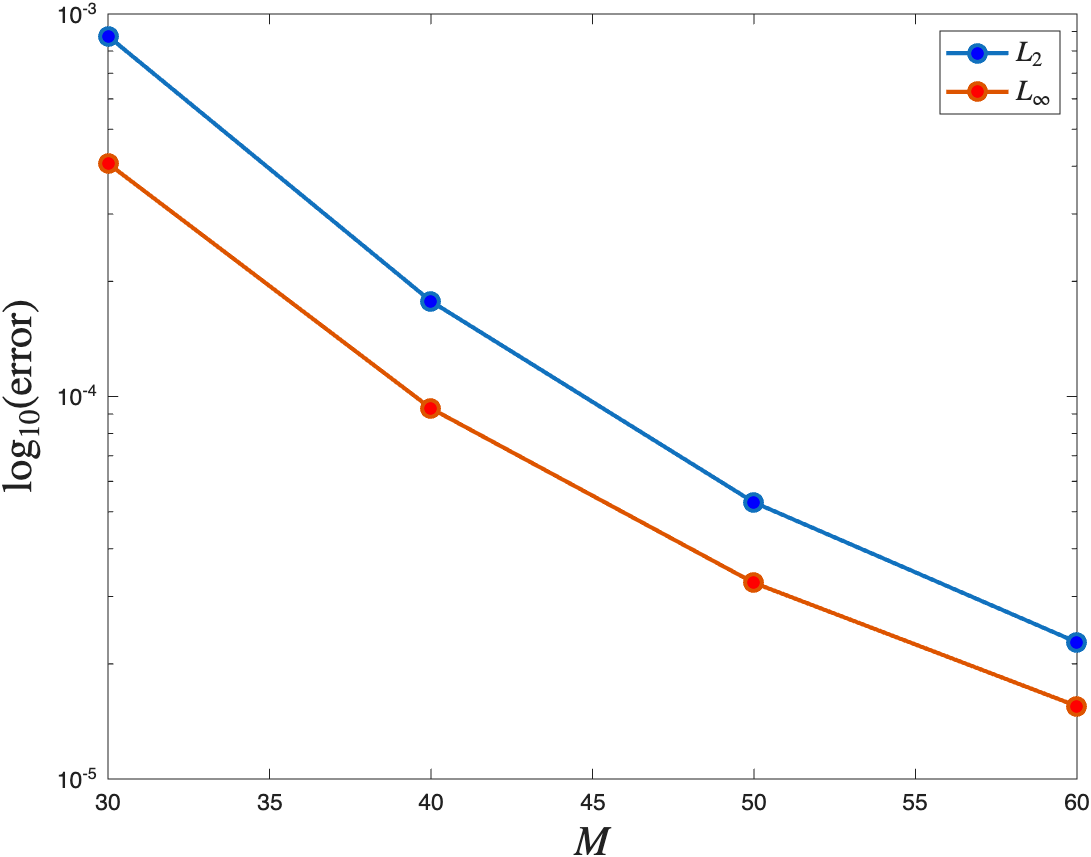}
        \caption{}
        \label{fig: ep2-error-x}
    \end{subfigure}
    \begin{subfigure}{0.32\textwidth}
        \centering
        \includegraphics[width=\textwidth]{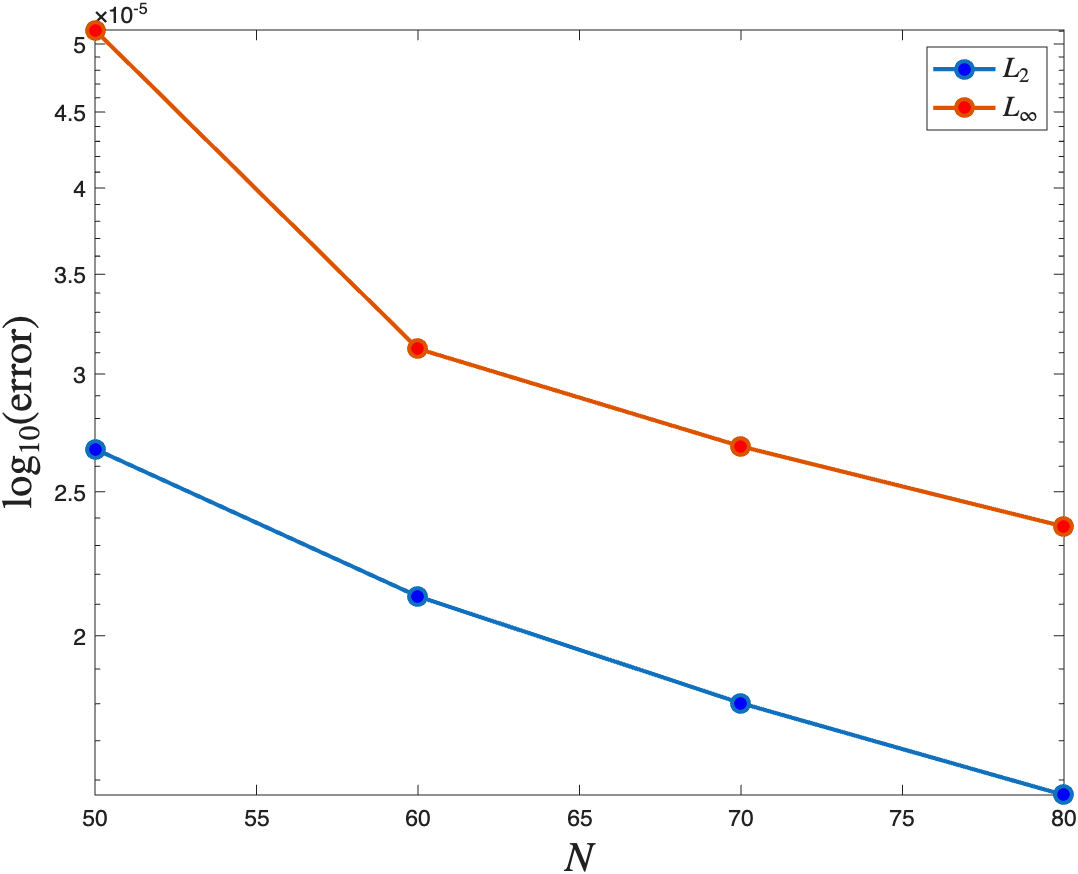}
        \caption{}
        \label{fig: ep2-error-xi}
    \end{subfigure}
    \begin{subfigure}{0.32\textwidth}
        \centering
        \includegraphics[width=\textwidth]{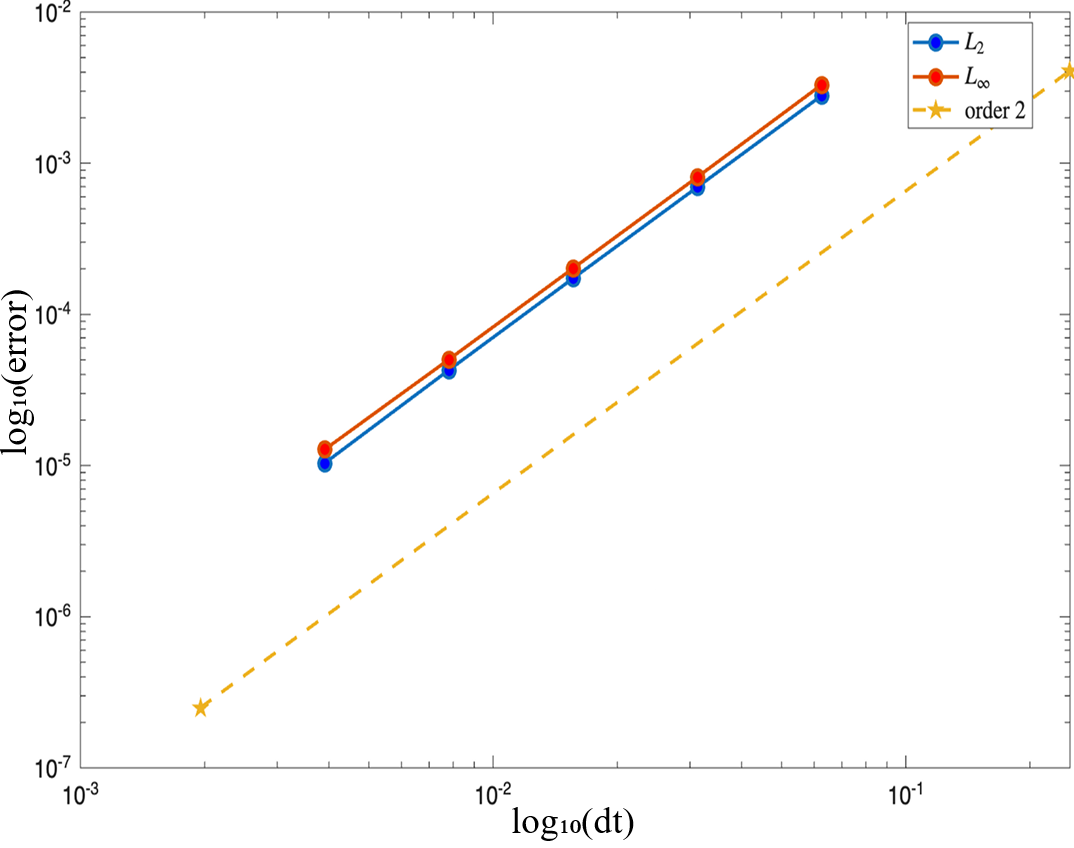}
        \caption{}
        \label{fig: ep2-error-time}
    \end{subfigure}
    \caption{Convergence rates for \textbf{Example 2} at $T=0.5$ with respect to $M$ (left), $N$ (middle), and time step $\Delta t$ (right). The red and blue lines correspond to the $L^2$ and $L^{\infty}$ errors, respectively.}
    \label{fig: example 2-nonharmonic potential error order}
\end{figure}

As shown in Figure \ref{fig: example 2-nonharmonic potential error order}, the proposed TSSP method again exhibits spectral convergence with respect to $M$ and $N$, and second-order convergence with respect to the time step $\Delta t$.

\textbf{Example 3: WPFP Equation with $\alpha=-1$}

Now, we test the convergence of the proposed TSSP method for the WPFP system with parameter $\alpha=-1$. The computational domain is set to $[-4,4]\times[-4,4]$, and the total simulation time is $T=0.25$. All other parameters and the initial conditions are the same as in Example 1. The reference solution was computed on a $2^{10} \times 2^{10}$ spatial grid with a time step of $\Delta t=2^{-10}$.

\begin{figure}[htbp]
    \centering
    \begin{subfigure}{0.32\textwidth}
        \centering
        \includegraphics[width=\textwidth]{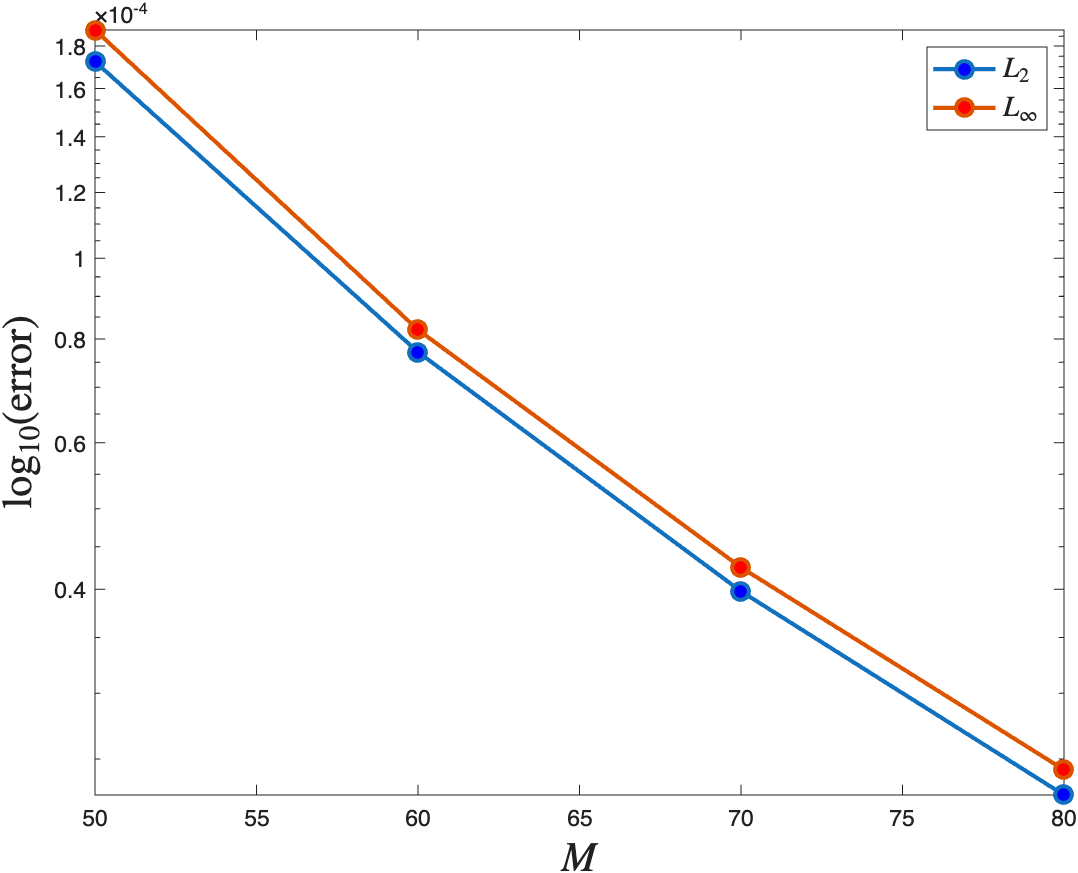}
        \caption{}
        \label{fig: ep3-error-x}
    \end{subfigure}
    \begin{subfigure}{0.32\textwidth}
        \centering
        \includegraphics[width=\textwidth]{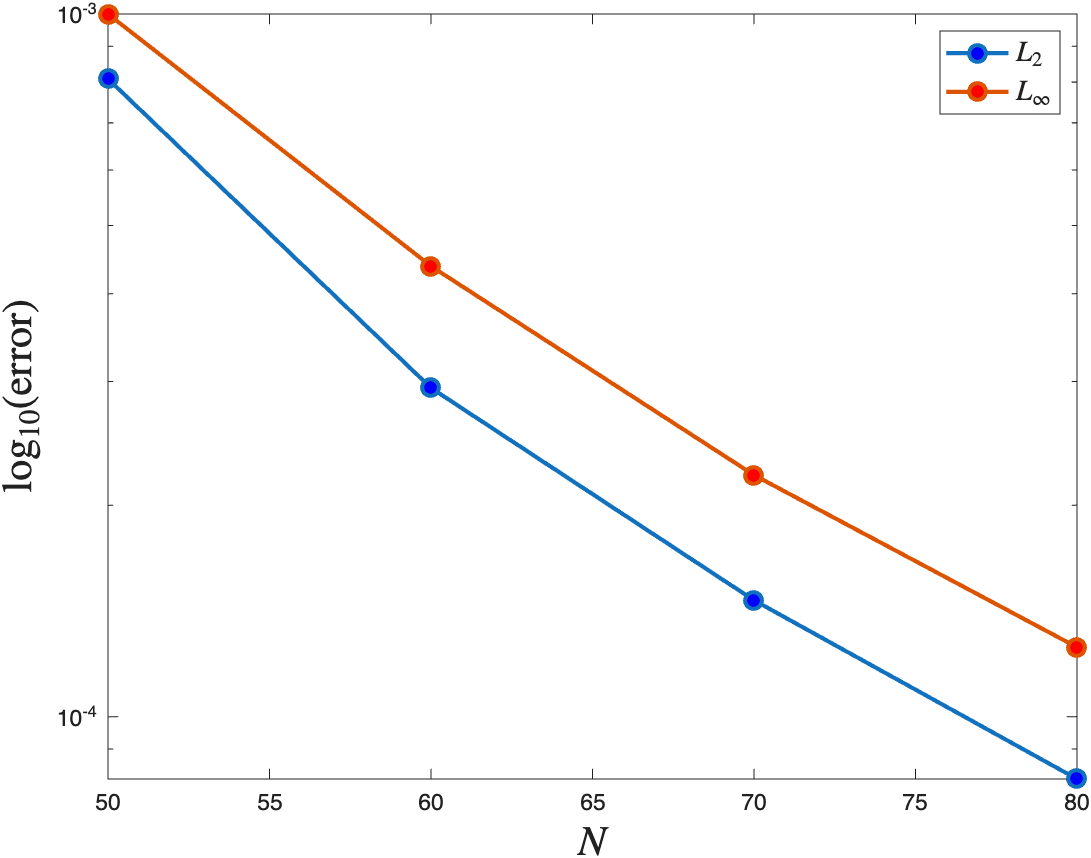}
        \caption{}
        \label{fig: ep3-error-xi}
    \end{subfigure}
    \begin{subfigure}{0.32\textwidth}
        \centering
        \includegraphics[width=\textwidth]{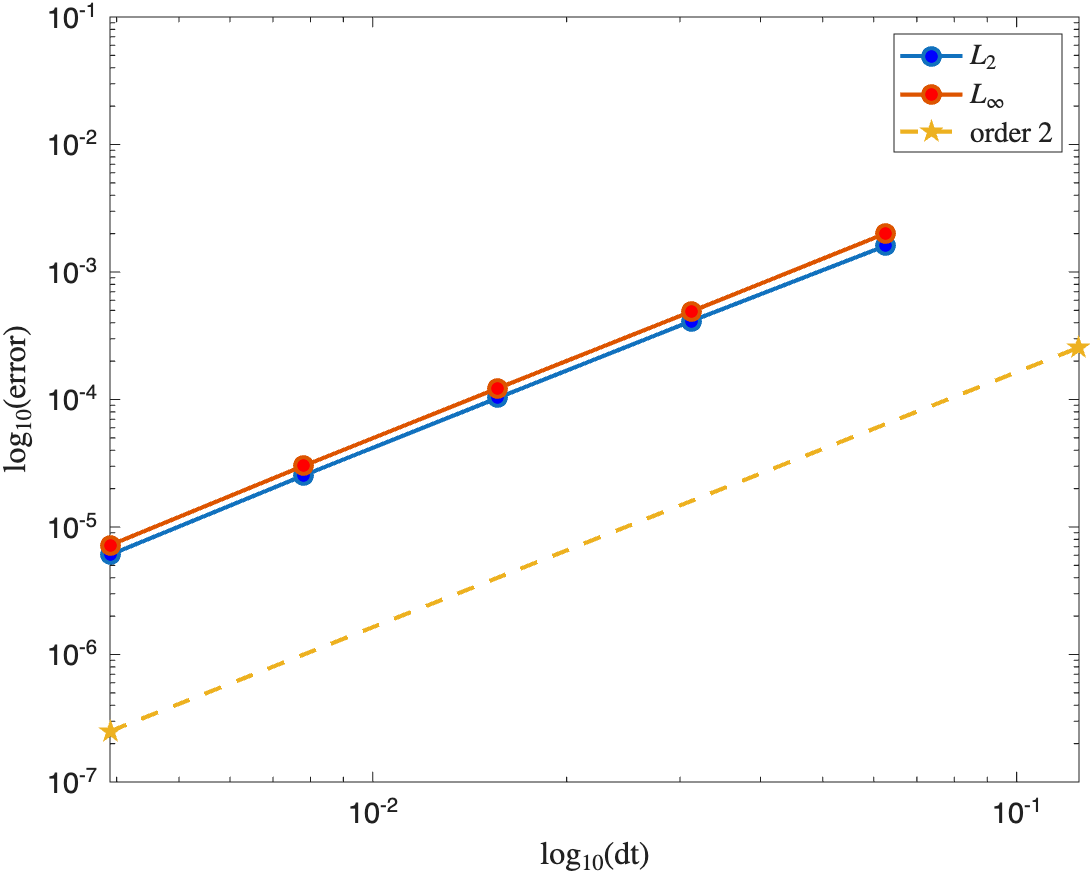}
        \caption{}
        \label{fig: ep3-error-time}
    \end{subfigure}
    \caption{Convergence rates for \textbf{Example 3} at $T=0.25$ with respect to $M$ (left), $N$ (middle), and time step $\Delta t$ (right). The red and blue lines correspond to the $L^2$ and $L^{\infty}$ errors, respectively.}
    \label{fig: example 3-WPFP error order}
\end{figure}

The error plots in Figure \ref{fig: example 3-WPFP error order} verify the spectral accuracy of the TSSP method in space and momentum, along with its second-order convergence in time.

\subsection{Existence of Steady States in WFP and WPFP Systems}

\textbf{Example 4: Steady state of the WFP equation}

This example investigates the evolution to a steady state of the WFP system for two types of potentials: one is adding a small perturbation to a harmonic potential (near-harmonic) and the other deviates significantly from harmonic potential (far-from-harmonic). The parameters are set to $\varepsilon=0.1$, $D_{pp}=D_{qq}=0.1$, $D_{pq}=0.0$, and $\gamma=1$. The initial conditions are $a_{11}=a_{22}=-1$, $x_0=0.1$, and $\xi_0=-0.2$. The simulation is performed on a $[-4,4]\times[-4,4]$ domain with a $2^7 \times 2^7$ grid and a time step of $\Delta t=2^{-8}$.

\begin{figure}[htbp]
    \centering
    \begin{subfigure}{0.32\textwidth}
        \centering
        \includegraphics[width=\textwidth]{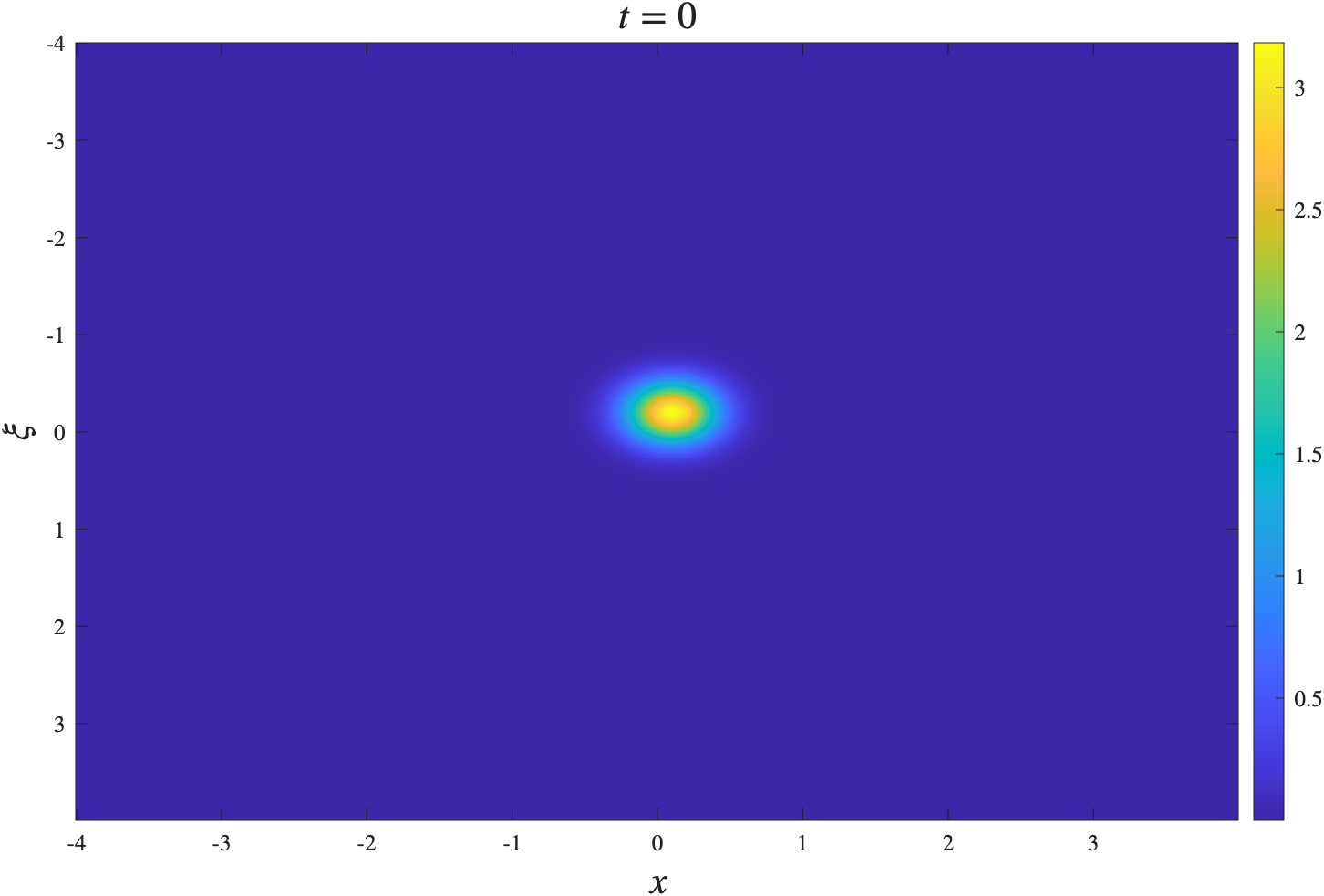}
        \label{fig: ep4-near-harmonic-Wigner-function-t0}
    \end{subfigure}
    \begin{subfigure}{0.32\textwidth}
        \centering
        \includegraphics[width=\textwidth]{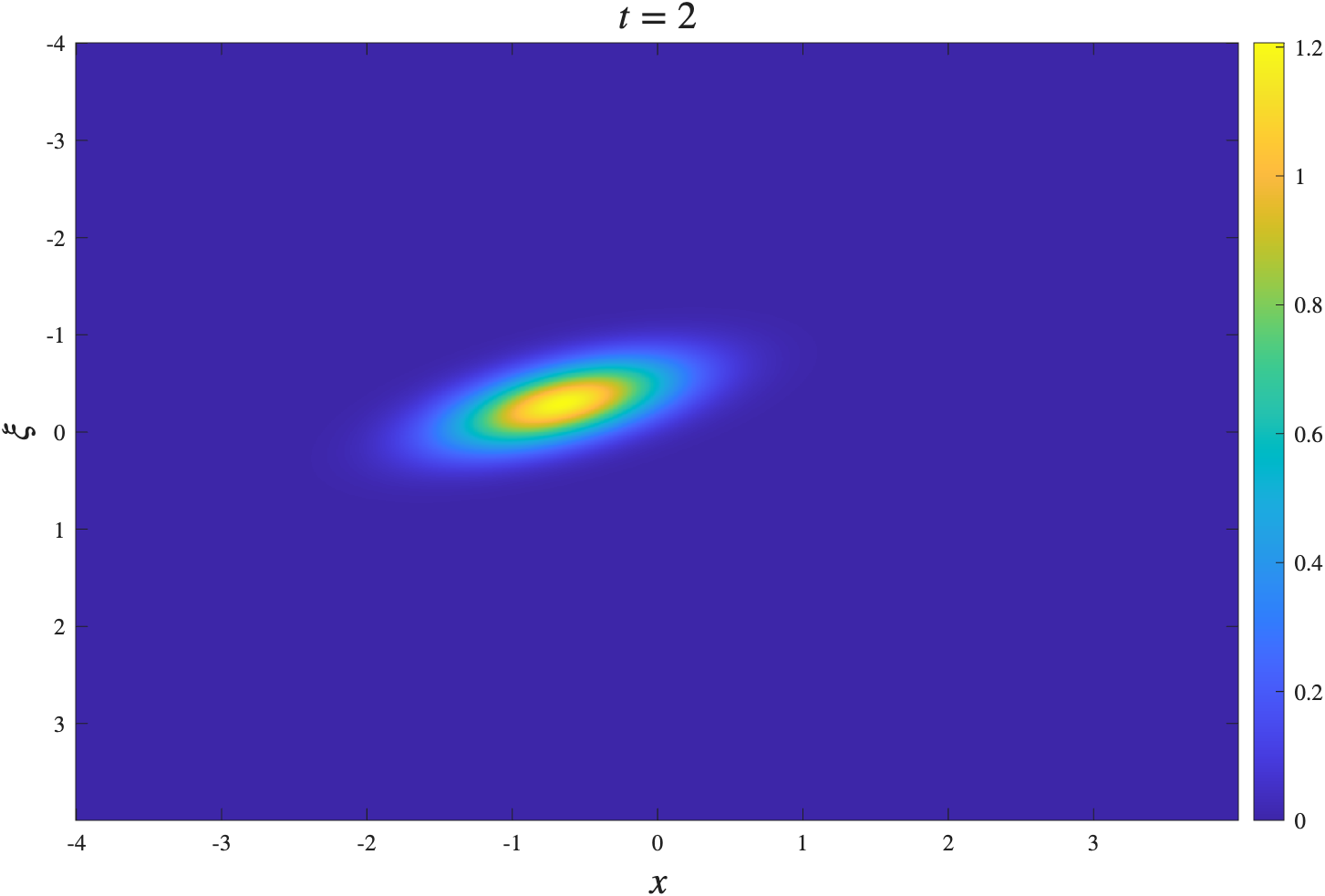}
        \label{fig: ep4-near-harmonic-Wigner-function-t2}
    \end{subfigure}
    \begin{subfigure}{0.32\textwidth}
        \centering
        \includegraphics[width=\textwidth]{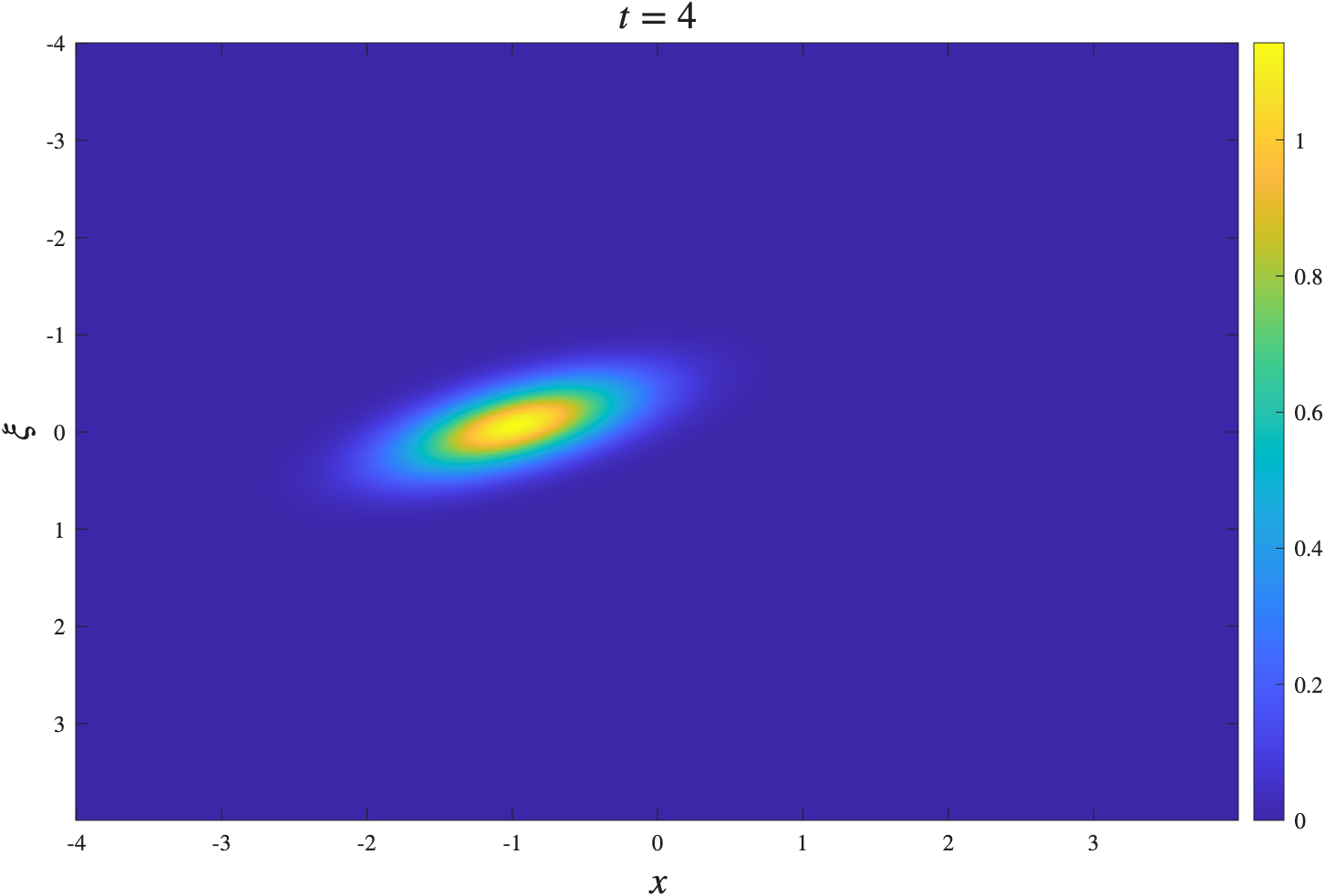}
        \label{fig: ep4-near-harmonic-Wigner-function-t4}
    \end{subfigure}
    \\
    \begin{subfigure}{0.32\textwidth}
        \centering
        \includegraphics[width=\textwidth]{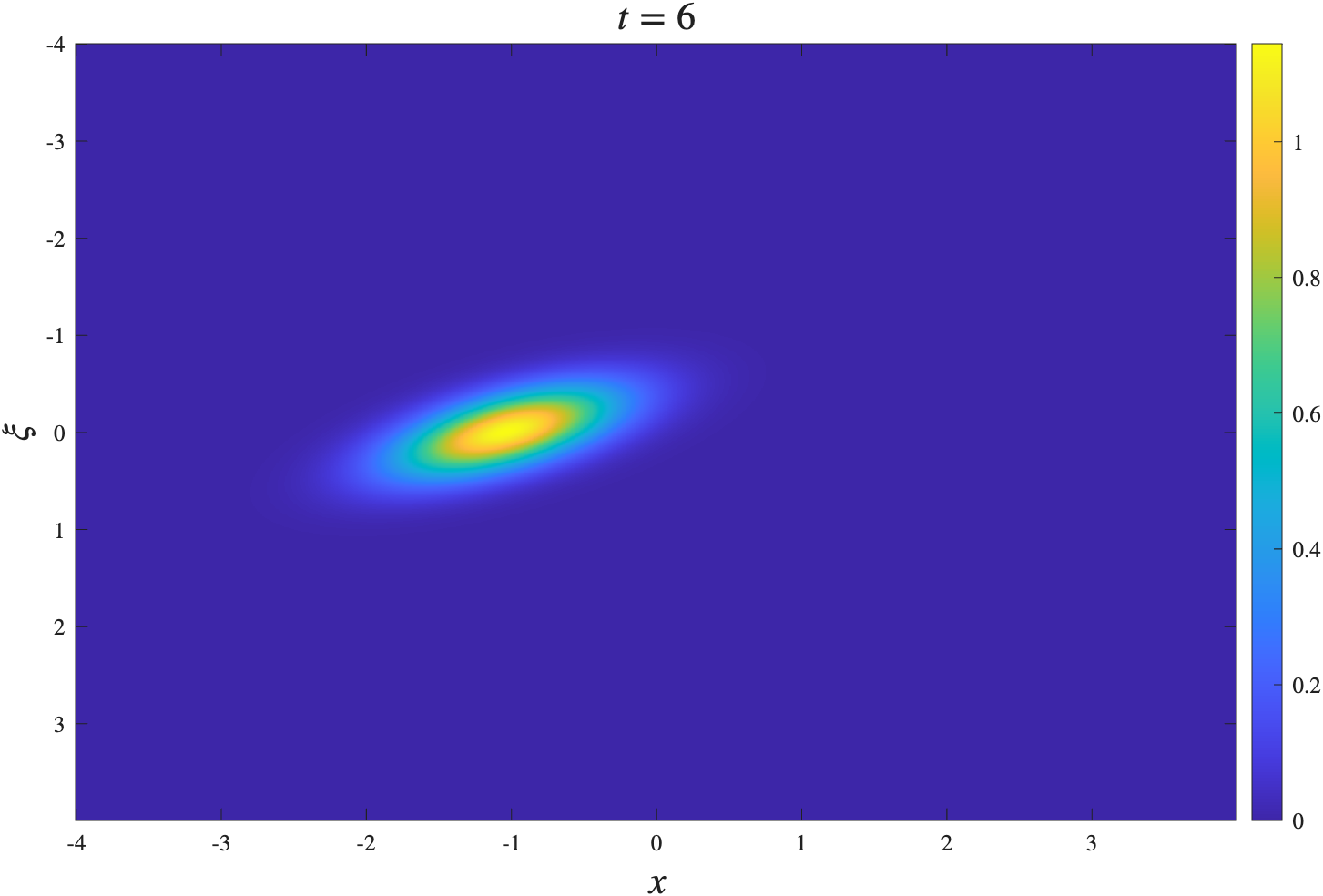}
        \label{fig: ep4-near-harmonic-Wigner-function-t6}
    \end{subfigure}
    \begin{subfigure}{0.32\textwidth}
        \centering
        \includegraphics[width=\textwidth]{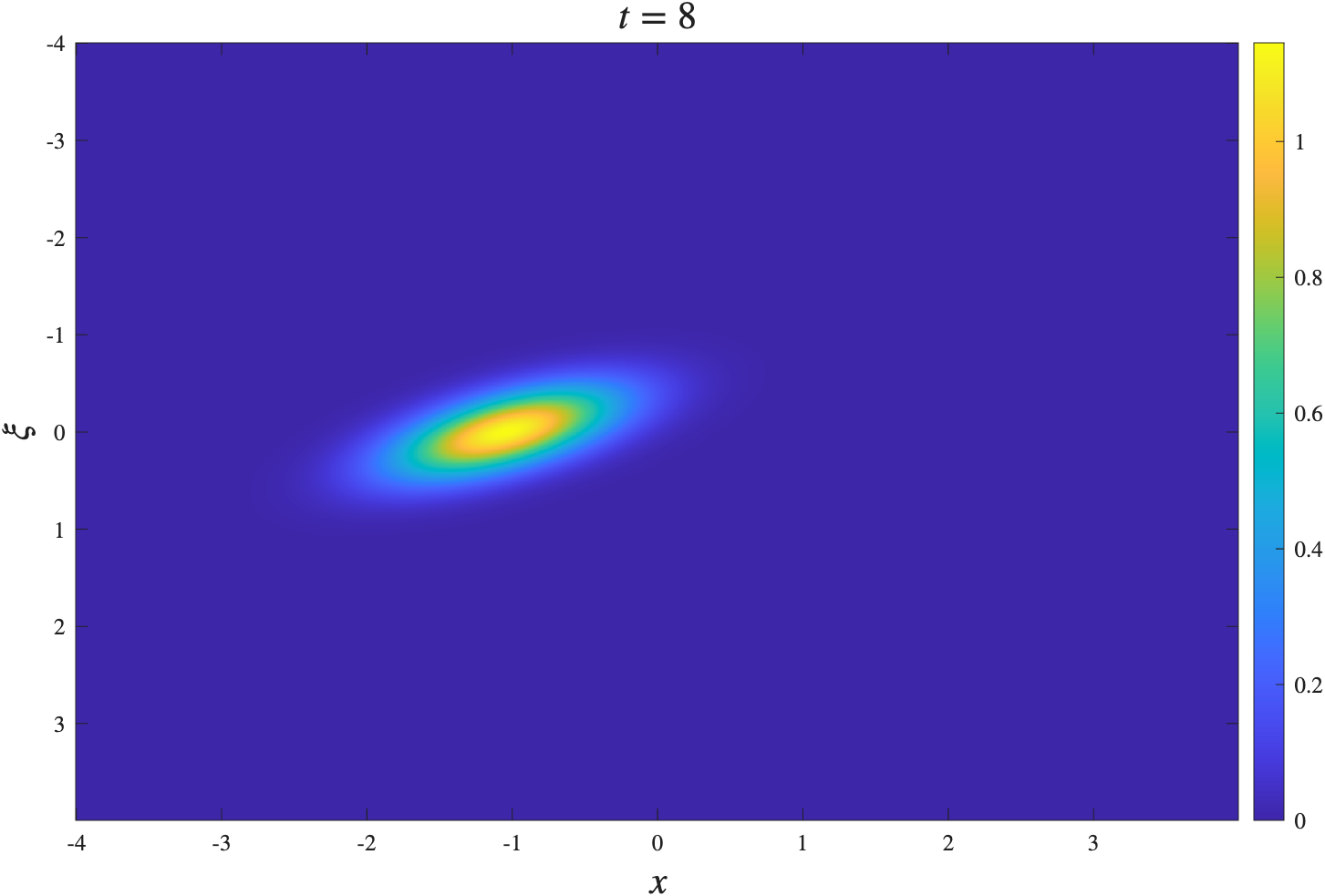}
        \label{fig: ep4-near-harmonic-Wigner-function-t8}
    \end{subfigure}
    \begin{subfigure}{0.32\textwidth}
        \centering
        \includegraphics[width=\textwidth]{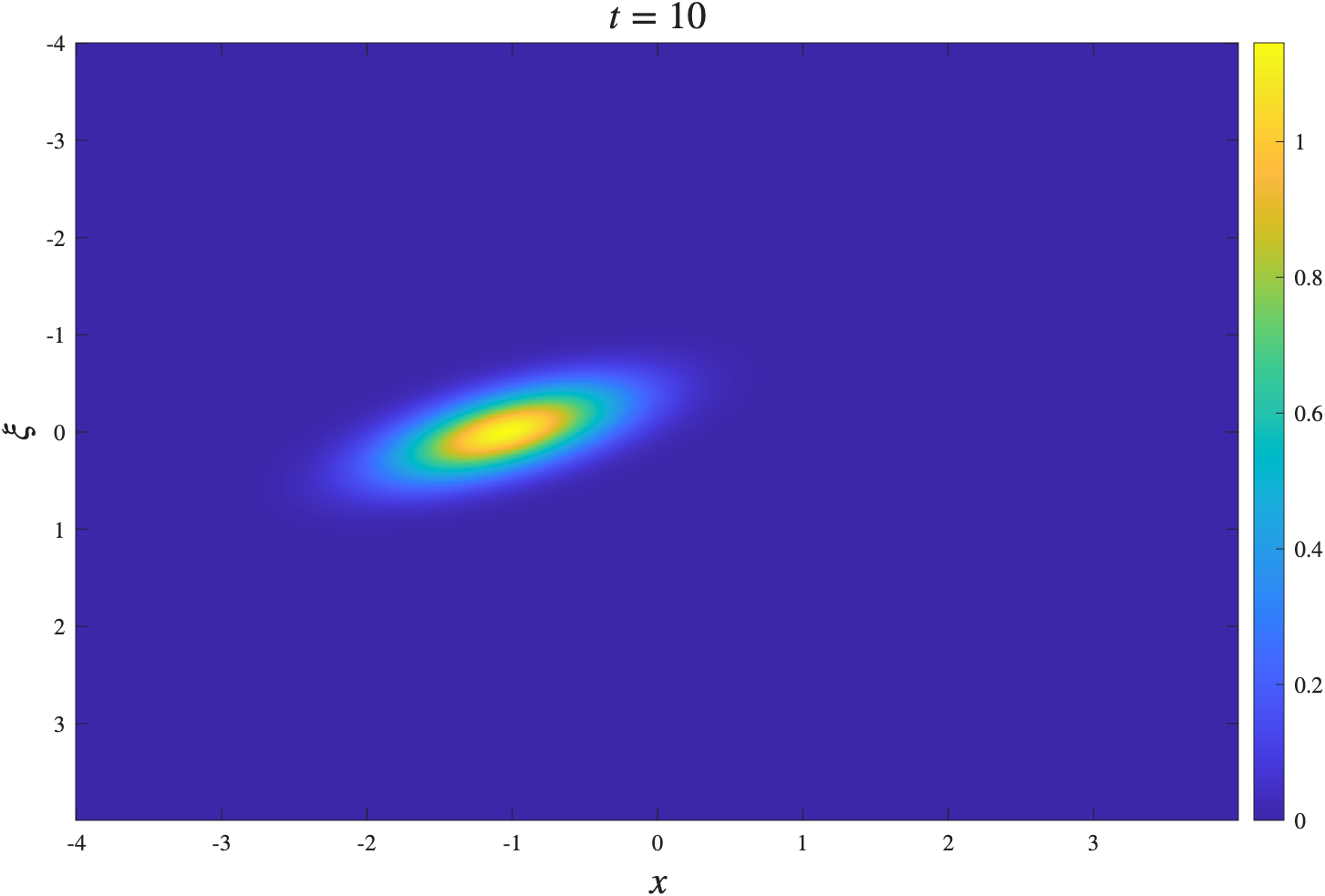}
        \label{fig: ep4-near-harmonic-Wigner-function-t10}
    \end{subfigure}
    \caption{Snapshots of the Wigner function evolution for \textbf{Example 4} with the near-harmonic potential $V(x)=\frac{1}{2}x^2+x+0.1\sin(x)$. The panels display the function at times $t=0,2,4,6,8$, and $10$.}
    \label{fig: ep4-near-harmonic-Wigner-function}
\end{figure}

\begin{figure}[htbp]
    \centering
    \includegraphics[width=\textwidth]{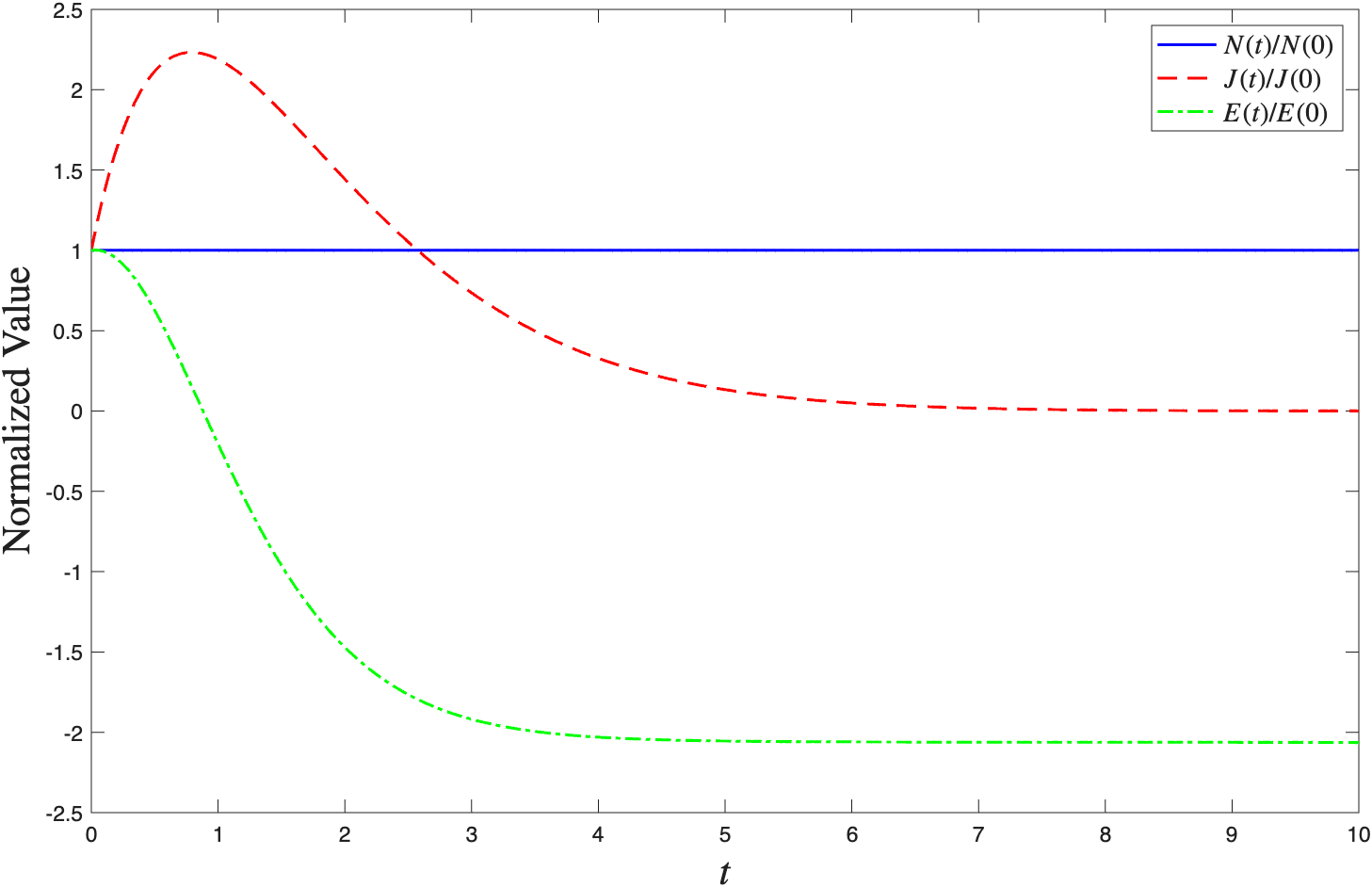}
    \caption{Evolution of normalized physical quantities for \textbf{Example 4} with the near-harmonic potential. The total particle number $N(t)$, total momentum $J(t)$, and total energy $E(t)$ are normalized by their initial values.}
    \label{fig: ep4-near-harmonic-physical-quantity}
\end{figure}

First, we consider a potential near the harmonic form:
\[ V(x)=\frac{1}{2}x^2+x+0.1\sin(x). \]
As depicted in Figure \ref{fig: ep4-near-harmonic-Wigner-function}, the Wigner function evolves and stabilizes after $t\rightarrow8$. This is further supported by Figure \ref{fig: ep4-near-harmonic-physical-quantity}, where physical quantities $N(t), J(t)$ and $E(t)$  plateau over time. These findings are consistent with the theoretical results presented in \cite{arnold2004analysis}.

\begin{figure}[htbp]
    \centering
    \begin{subfigure}{0.48\textwidth}
        \centering
        \includegraphics[width=\textwidth]{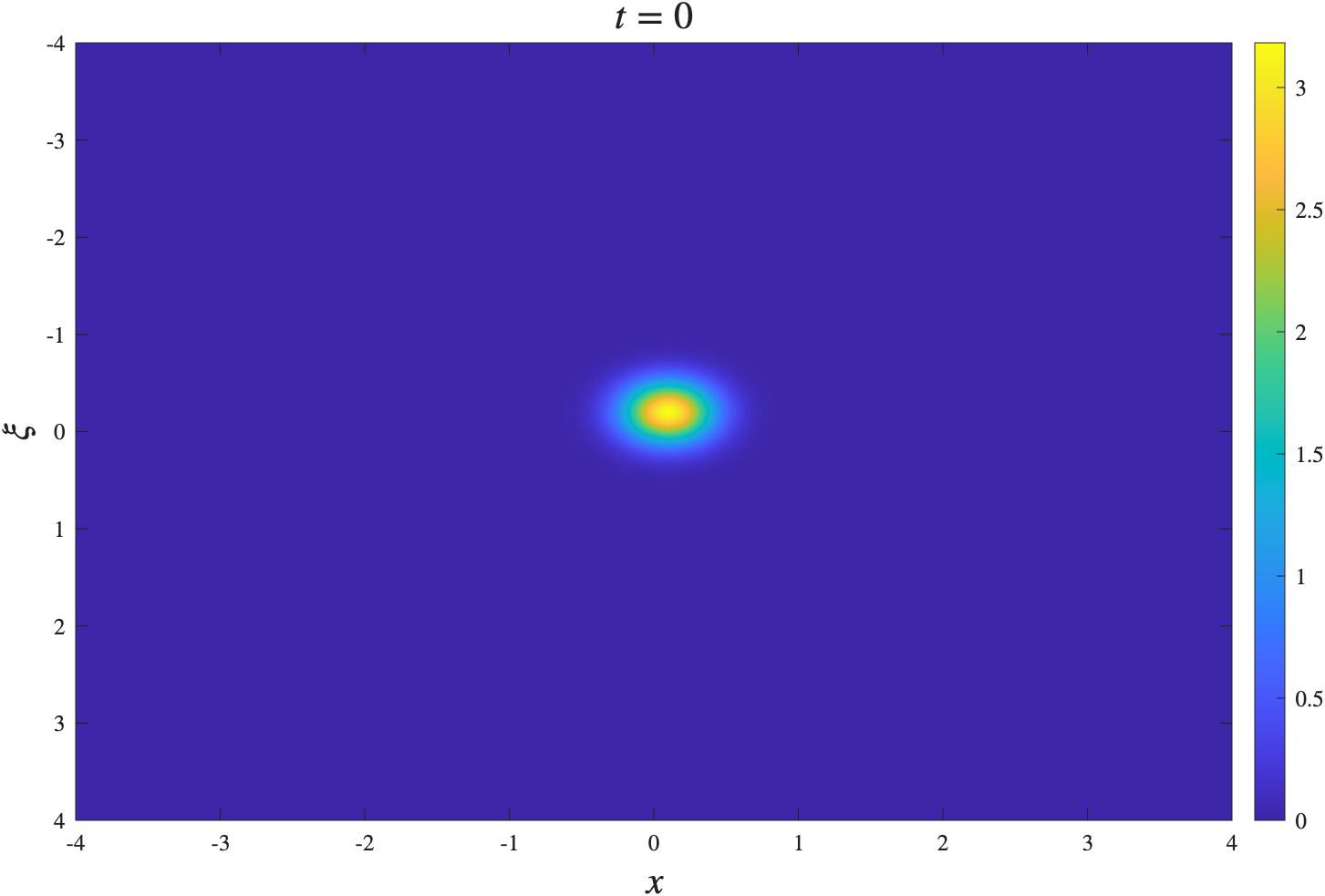}
        \label{fig: ep4-far-harmonic-Wigner-function-t0}
    \end{subfigure}
    \begin{subfigure}{0.48\textwidth}
        \centering
        \includegraphics[width=\textwidth]{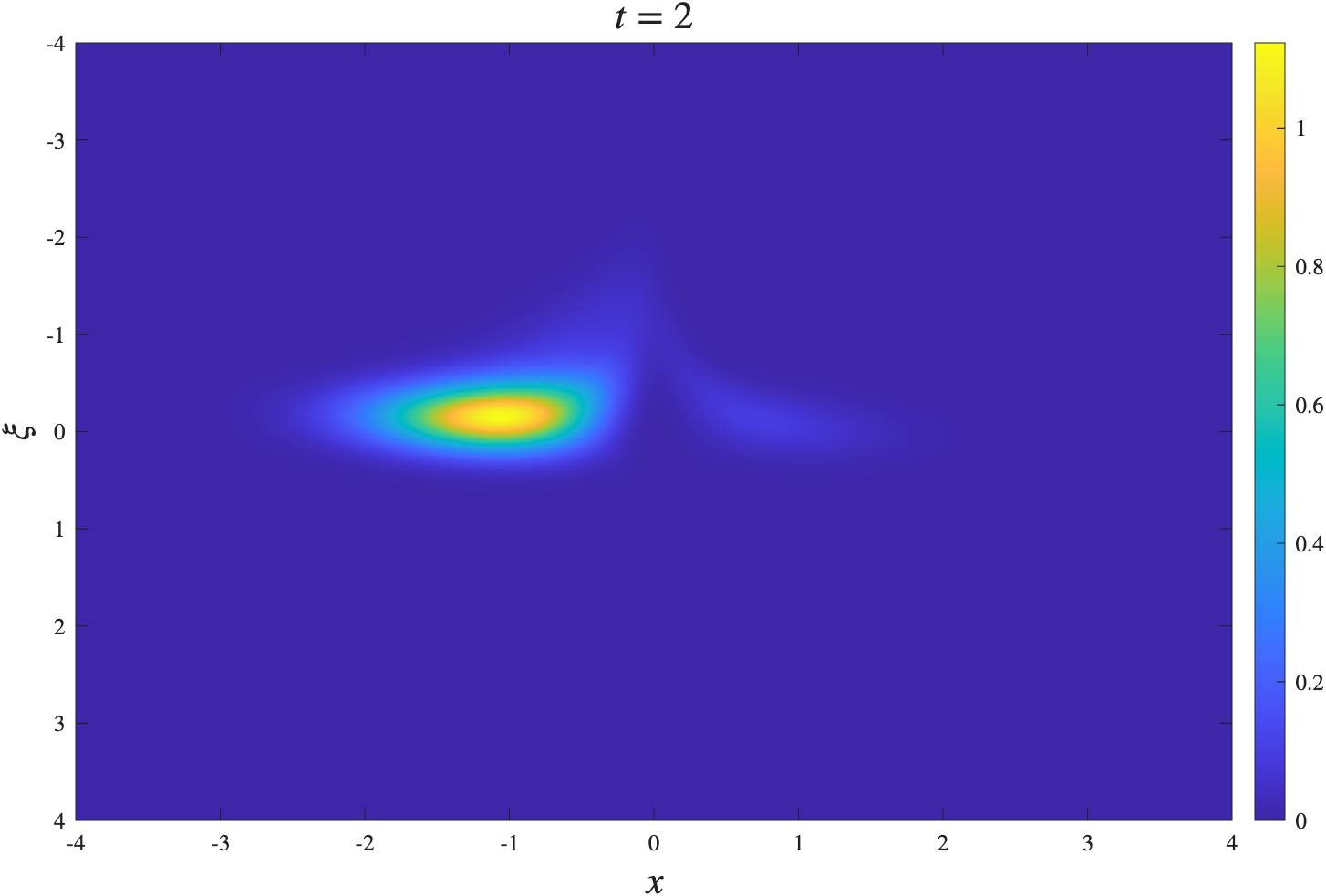}
        \label{fig: ep4-far-harmonic-Wigner-function-t2}
    \end{subfigure}
    \\
    \begin{subfigure}{0.48\textwidth}
        \centering
        \includegraphics[width=\textwidth]{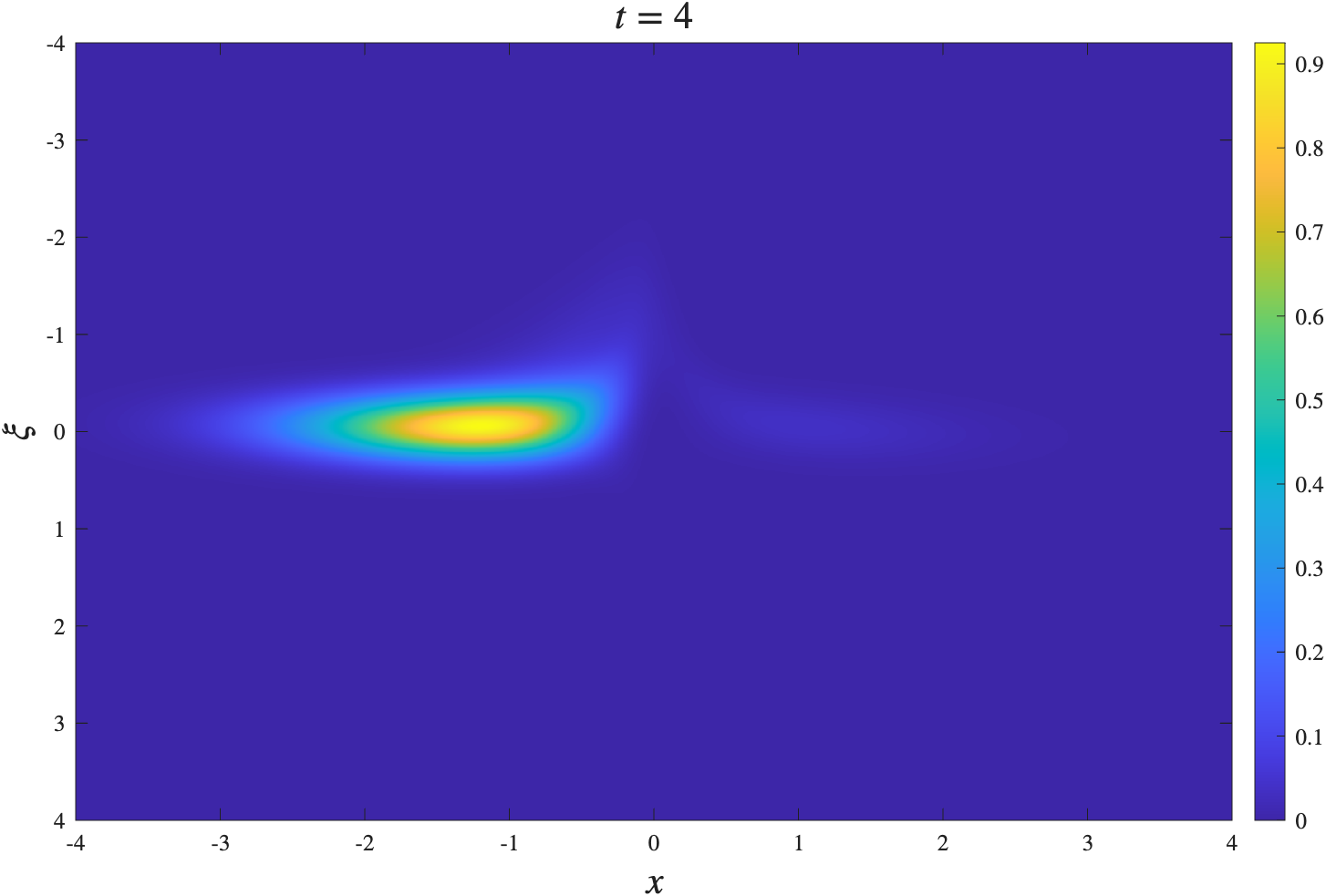}
        \label{fig: ep4-far-harmonic-Wigner-function-t4}
    \end{subfigure}
    \begin{subfigure}{0.48\textwidth}
        \centering
        \includegraphics[width=\textwidth]{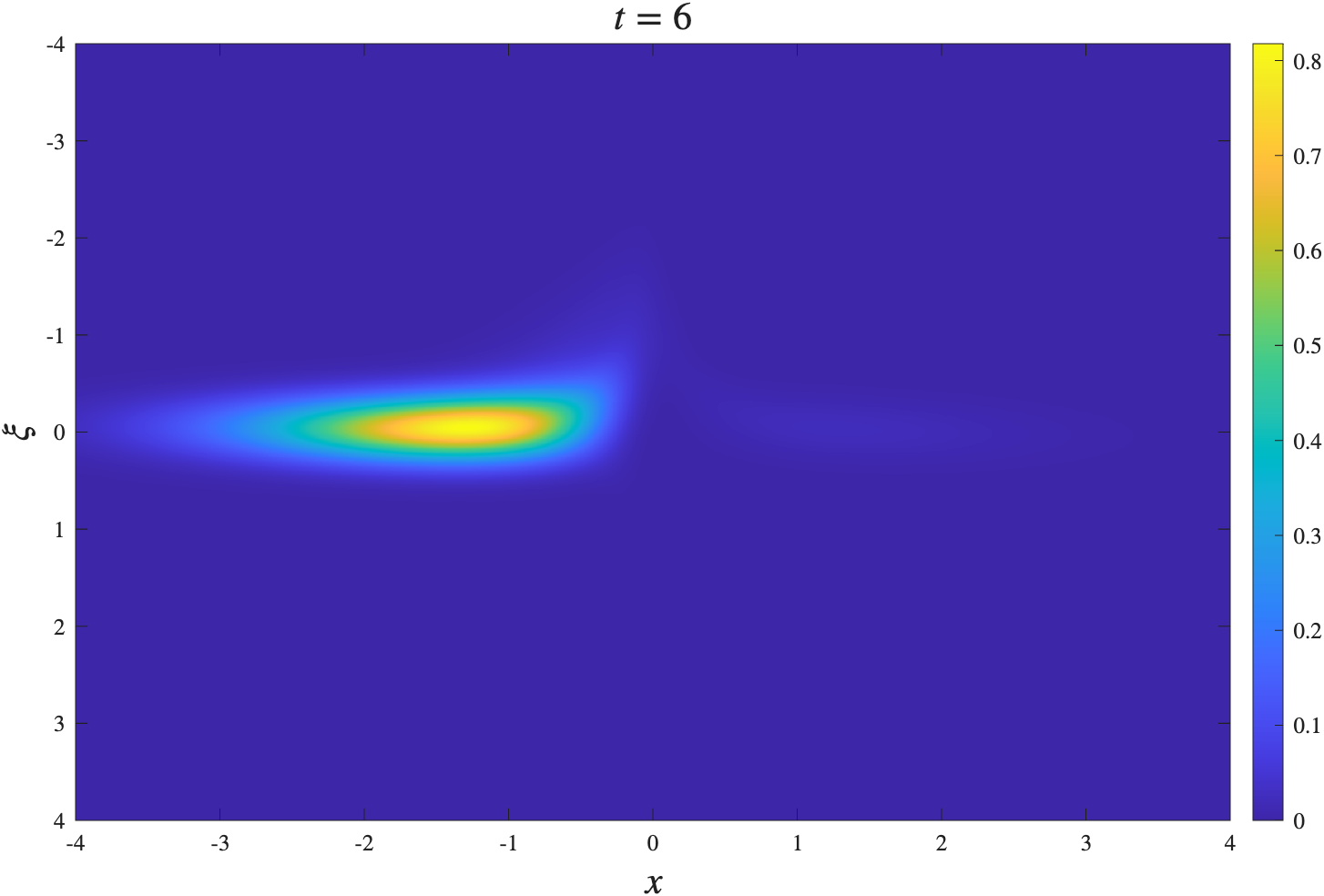}
        \label{fig: ep4-far-harmonic-Wigner-function-t6}
    \end{subfigure}
    \caption{Snapshots of the Wigner function evolution for \textbf{Example 4} with the far-from-harmonic potential $V(x)=\arctan(10x)+\frac{\pi}{2}$. The panels display the function at times $t=0,2,4$, and $6$.}
    \label{fig: ep4-far-harmonic-Wigner-function}
\end{figure}

\begin{figure}[htbp]
    \centering
    \includegraphics[width=\textwidth]{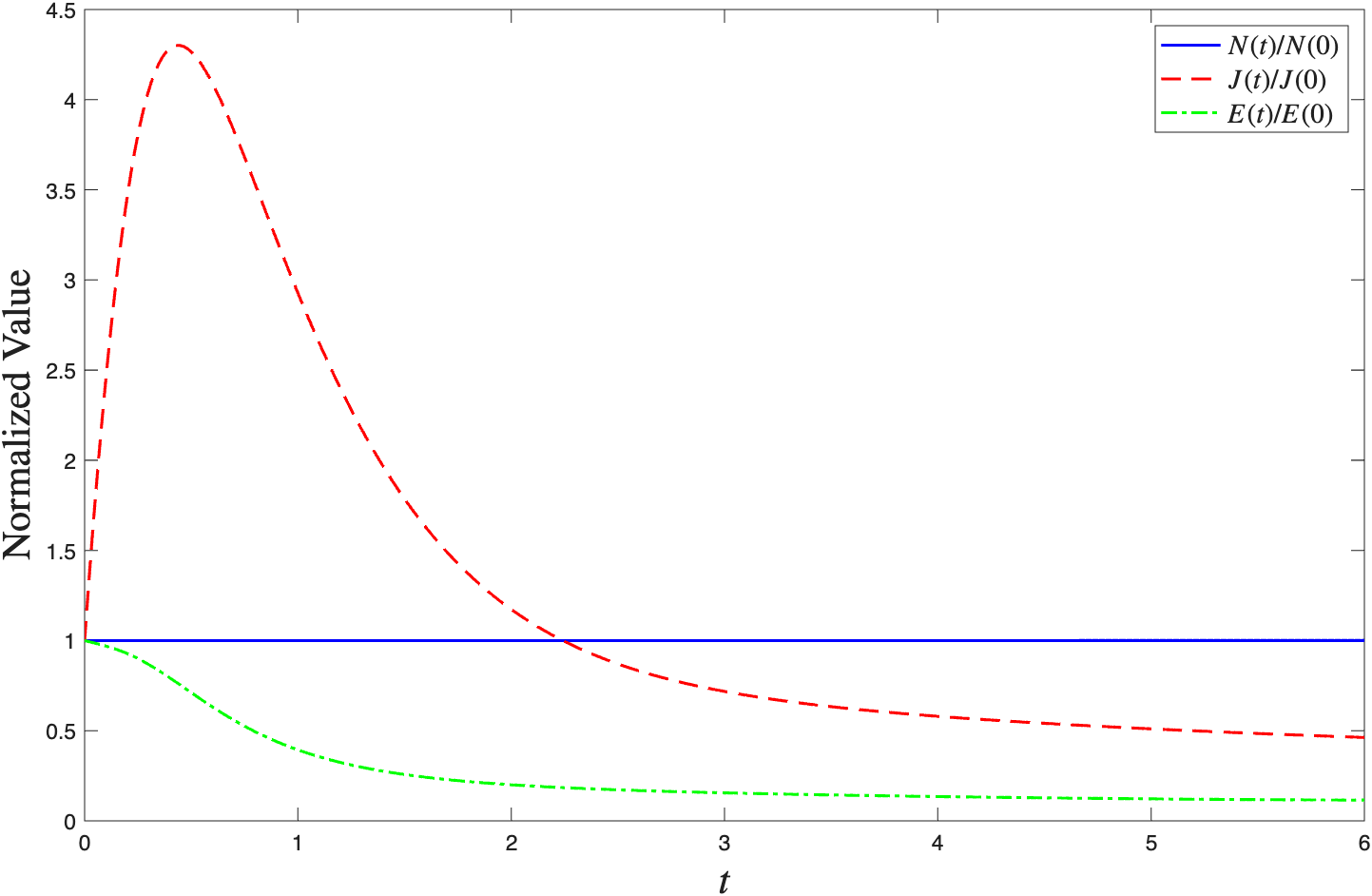}
    \caption{Evolution of normalized physical quantities for \textbf{Example 4} with the far-from-harmonic potential. The total particle number $N(t)$, total momentum $J(t)$, and total energy $E(t)$ are normalized by their initial values.}
    \label{fig: ep4-far-harmonic-physical-quantity}
\end{figure}

Next, we choose a potential far from the harmonic form:
\[ V(x)=\arctan(10x)+\frac{\pi}{2}. \]
Figures \ref{fig: ep4-far-harmonic-Wigner-function} and \ref{fig: ep4-far-harmonic-physical-quantity} present the evolution of the Wigner function and the corresponding physical quantities, respectively. From these figures, we observe that the system again reaches a steady state, with the Wigner function and the important physical quantities becoming constant after $t \approx 4$. This numerical result suggests that a steady state can be achieved even for potentials far from the harmonic form. This finding, to our best knowledge, has not yet been theoretically established in the literature.

\textbf{Example 5: Steady State of the WPFP Equation}

Finally, we numerically investigate the existence of the steady state for the WPFP equation. The model parameters are set as $\varepsilon=1$, $D_{pp}=D_{qq}=0.3$, $D_{pq}=0.0$, and $\gamma=1$. The initial conditions are $a_{11}=a_{22}=-1$, $x_0=0.1$, and $\xi_0=0.1$. The simulation is performed on a $[-20,20]\times[-20,20]$ domain with a $2^7 \times 2^7$ grid and a time step of $\Delta t=2^{-8}$.

\begin{figure}[htbp]
    \centering
    \begin{subfigure}{0.32\textwidth}
        \centering
        \includegraphics[width=\textwidth]{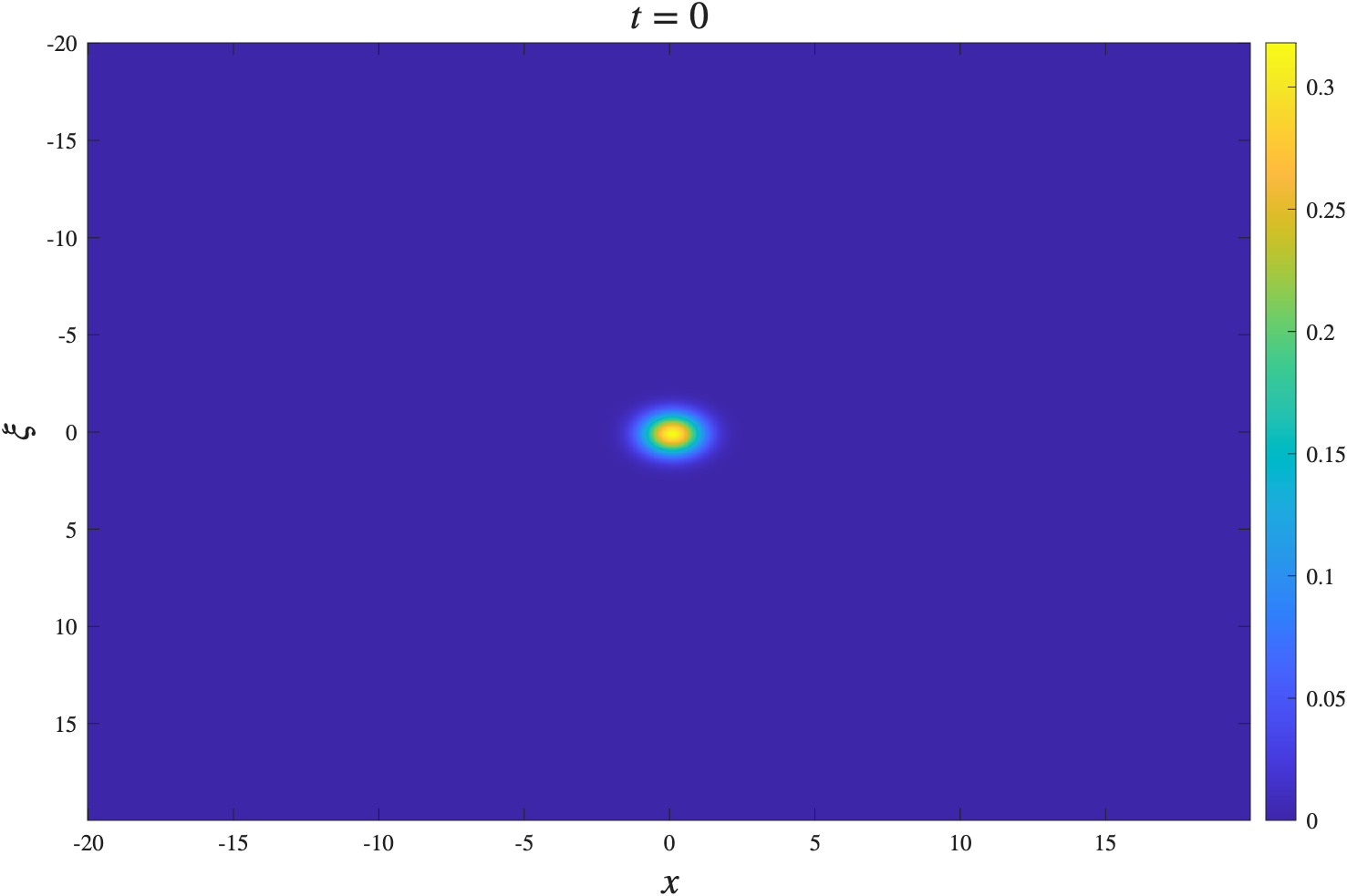}
        \label{fig: ep5-Wigner-function-t0}
    \end{subfigure}
    \begin{subfigure}{0.32\textwidth}
        \centering
        \includegraphics[width=\textwidth]{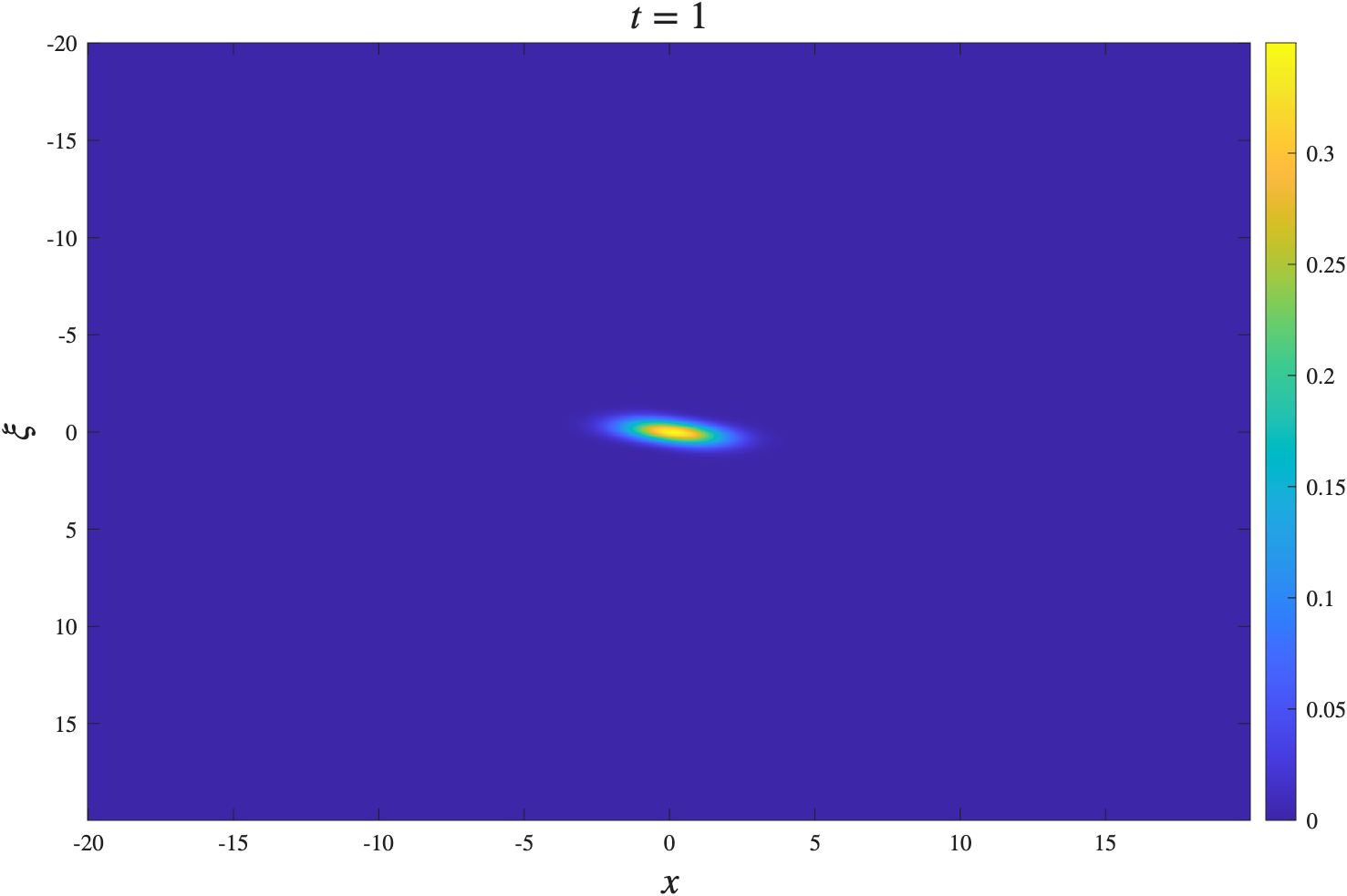}
        \label{fig: ep5-Wigner-function-t1}
    \end{subfigure}
    \begin{subfigure}{0.32\textwidth}
        \centering
        \includegraphics[width=\textwidth]{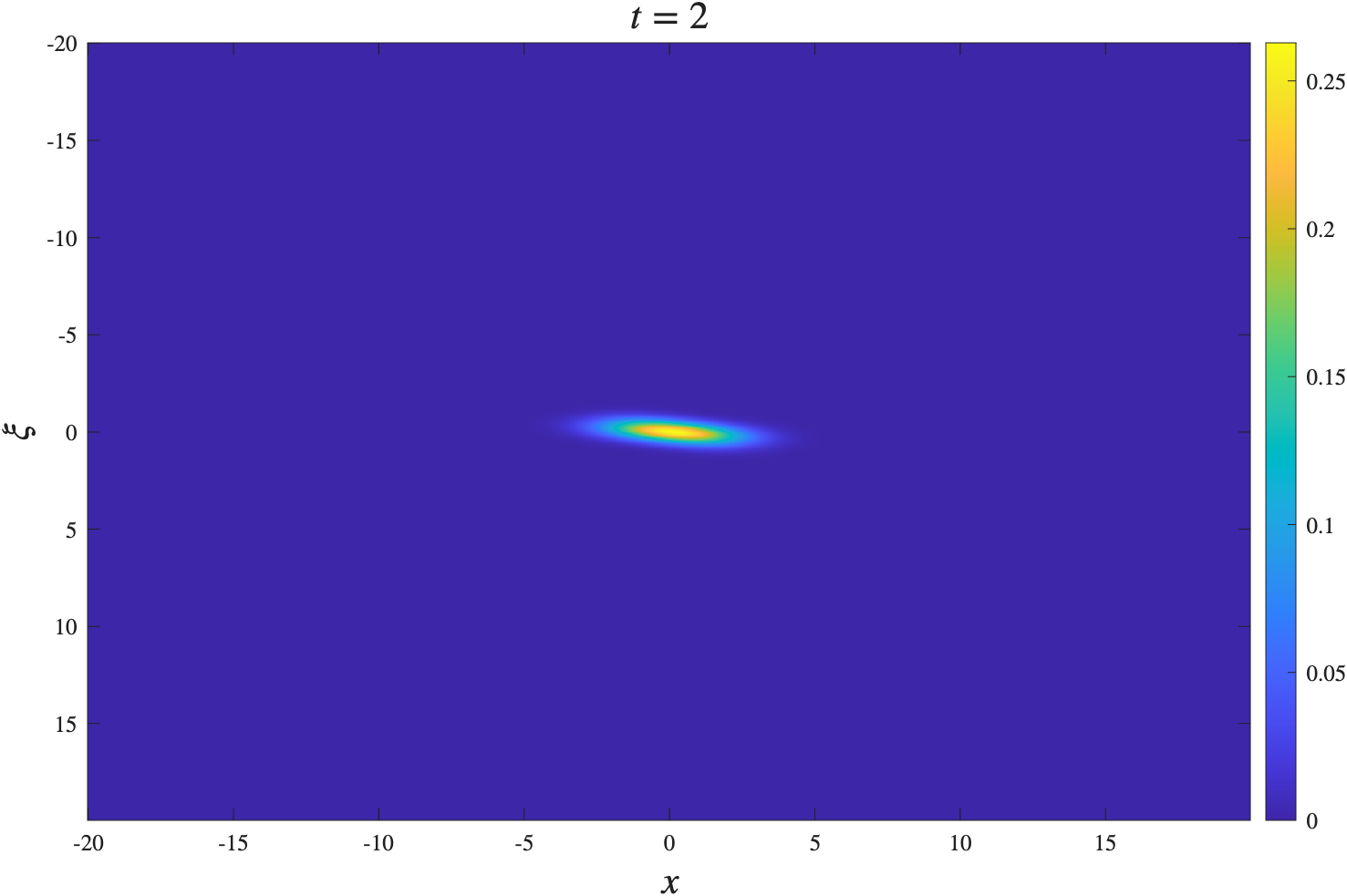}
        \label{fig: ep5-Wigner-function-t2}
    \end{subfigure}
    \\
    \begin{subfigure}{0.32\textwidth}
        \centering
        \includegraphics[width=\textwidth]{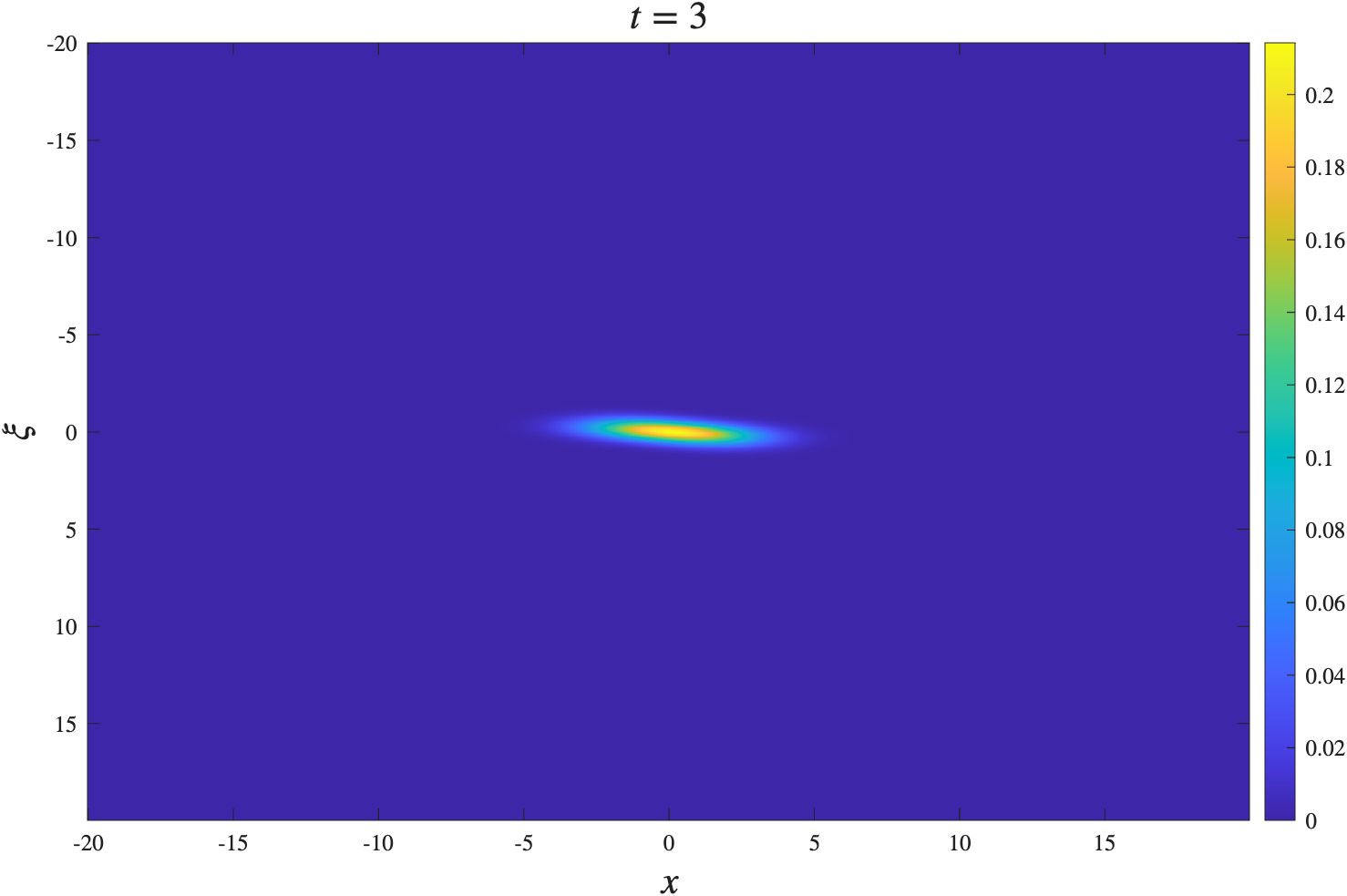}
        \label{fig: ep5-Wigner-function-t3}
    \end{subfigure}
    \begin{subfigure}{0.32\textwidth}
        \centering
        \includegraphics[width=\textwidth]{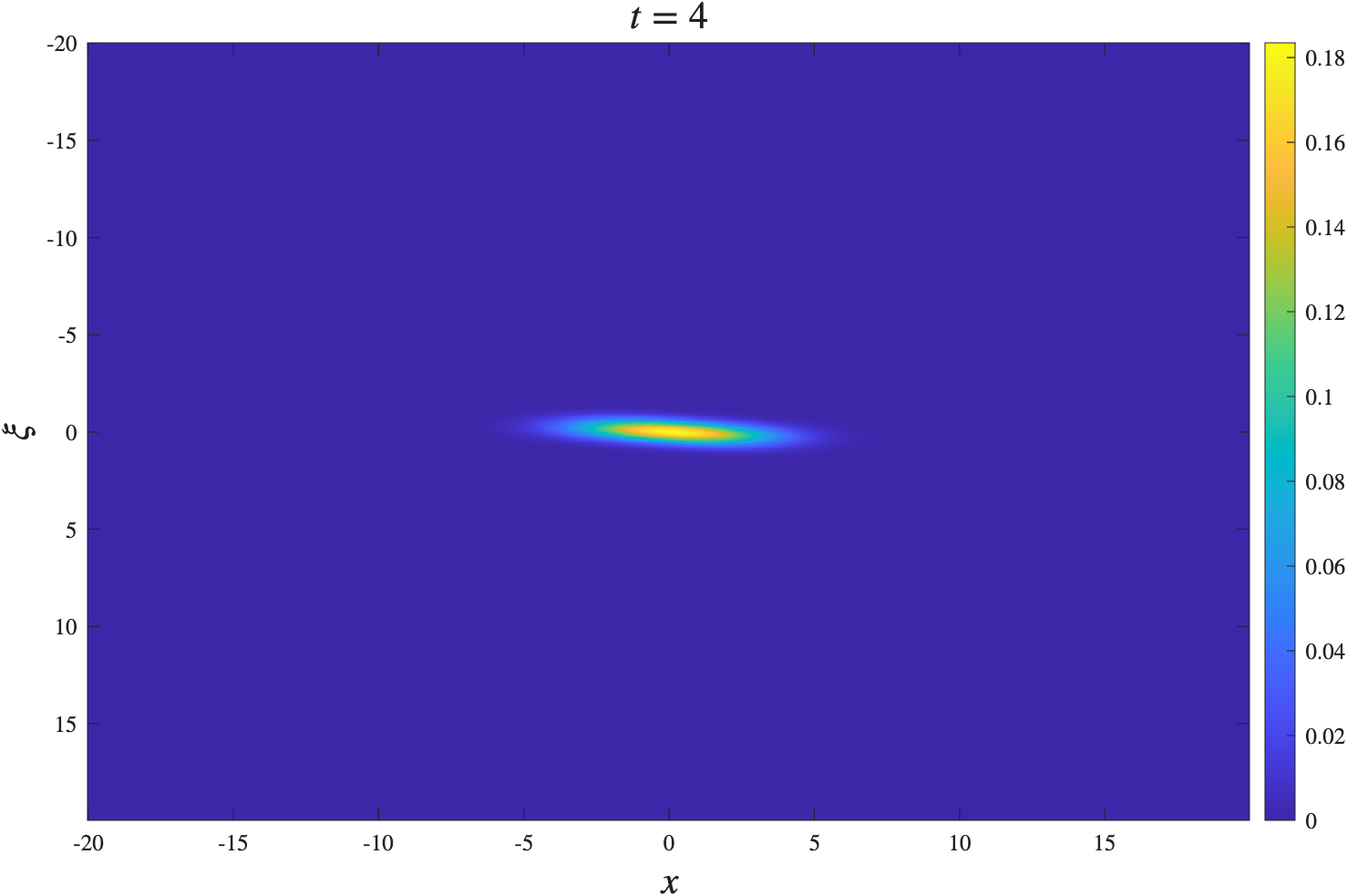}
        \label{fig: ep5-Wigner-function-t4}
    \end{subfigure}
    \begin{subfigure}{0.32\textwidth}
        \centering
        \includegraphics[width=\textwidth]{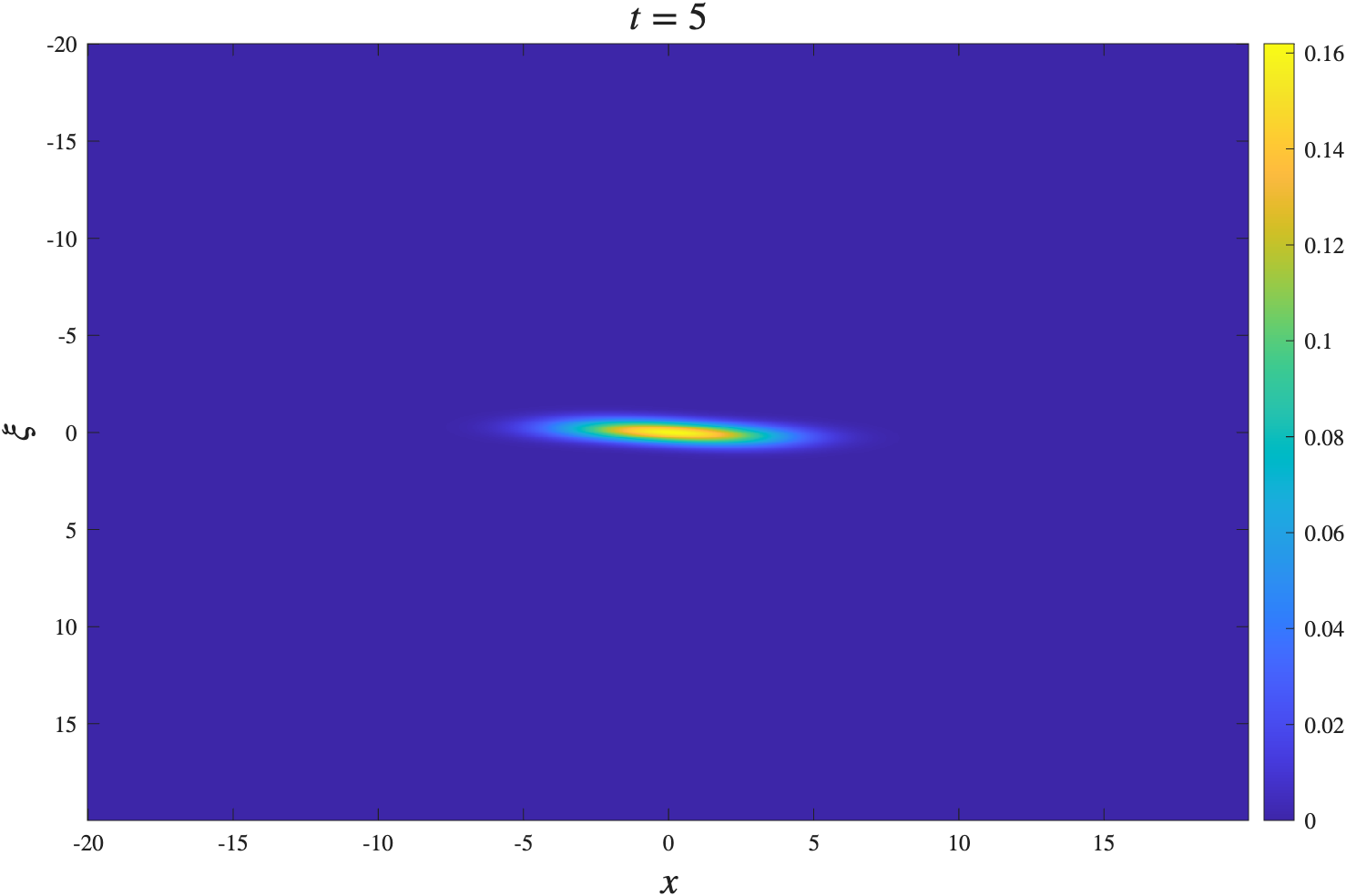}
        \label{fig: ep5-Wigner-function-t5}
    \end{subfigure}
    \caption{Snapshots of the Wigner function evolution for \textbf{Example 5}. The panels display the function at times $t=0,1,2,3,4$, and $5$, arranged sequentially from top-left to bottom-right.}
    \label{fig: ep5-Wigner-function}
\end{figure}

\begin{figure}[htbp]
    \centering
    \includegraphics[width=\textwidth]{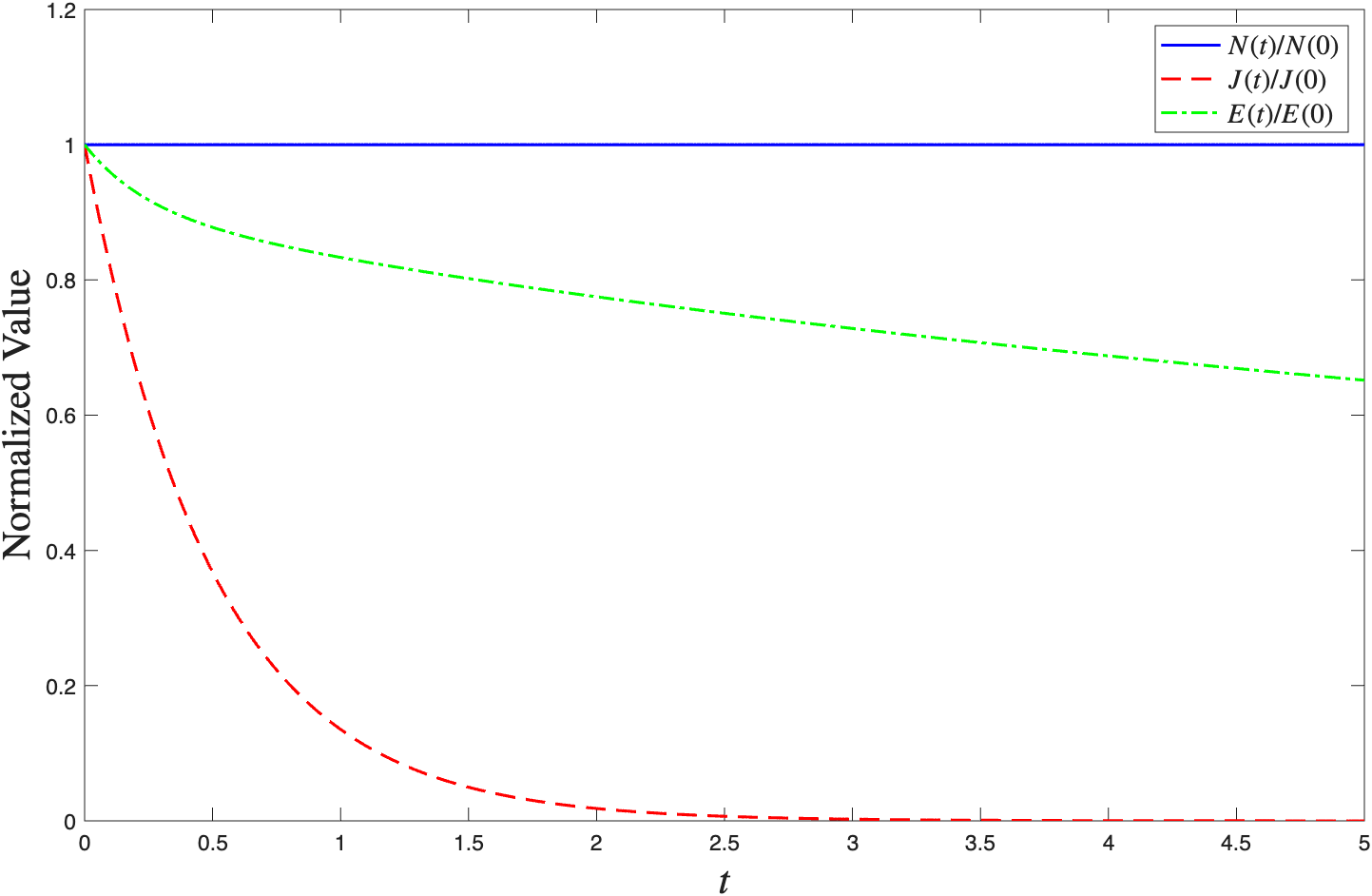}
    \caption{Evolution of normalized physical quantities for \textbf{Example 5}. The total particle number $N(t)$, total momentum $J(t)$, and total energy $E(t)$ are each normalized by their initial values.}
    \label{fig: ep5-physical-quantity}
\end{figure}

Figure \ref{fig: ep5-Wigner-function} displays snapshots of the Wigner function's evolution. It is observed that the function's profile ceases to change obviously after $t \approx 4$, indicating that the system has reached a steady state. Figure \ref{fig: ep5-physical-quantity} further corroborates this conclusion by showing that all the key physical quantities—total particle number, momentum, and energy—asymptotically approach constant values over time.

\section{Conclusion}
\label{sec: conclusion}

In this paper, we develop an efficient numerical scheme for the Wigner-Fokker-Planck and Wigner-Poisson-Fokker-Planck equations. The numerical method is constructed by combining a second-order time-splitting integrator with Fourier pseudospectral approximation, which leads to 
a scheme with second-order accuracy in time and spectral accuracy in phase space.
The accuracy and validation of the proposed method are verified by numerical tests.  The results reported by extensive numerical tests confirm these expected convergence rates, illustrating the robustness of the proposed algorithm. We then applied this scheme to simulate the long-time behavior of different quantum systems, where the convergence to steady states is successfully captured for both the WFP and the WPFP models. Notably, our results provide numerical evidence for the existence of steady states in WFP systems even with potentials far from the harmonic case, offering insight into phenomena unrevealed in theory.

This research opens several avenues for our future work. The most direct extension is the application of this method to higher-dimensional models, such as problems in $4$D phase space. To effectively mitigate the curse of dimensionality in such settings, we may combine the present scheme with some skillful methods dealing with high-dimensional problems, such as dynamical low-rank approximation methods. Another focus for future work will be put on the rigorous theoretical analysis of the developed numerical scheme.


\appendix

\section{Derivation of Fourier transform of nonlocal term}
\label{appendix: nonloacl}
\begin{align*}
&-\frac{1}{(\pi\varepsilon)^d}\int\int(\delta V)(x,y,t)\sum_{j=-\frac{N}{2}}^{\frac{N}{2}-1}\widehat{W}_j(x,t)\mathrm{e}^{\mathrm{i}\mu_j (\xi^{\prime}-c)}\mathrm{e}^{\frac{\mathrm{i}}{\varepsilon}2y\cdot\left(\xi^{\prime}-\xi\right)}\mathrm{d}\xi^{\prime}\mathrm{d}y\\
=&-\sum_{j=-\frac{N}{2}}^{\frac{N}{2}-1}\widehat{W}_j(x,t)\frac{1}{(\pi\xi)^d}\int\left[\int \mathrm{e}^{\mathrm{i}\xi^{\prime}(\mu_j+\frac{2y}{\varepsilon})}\mathrm{d}\xi^{\prime}\right](\delta V)(x,y,t)\mathrm{e}^{-\frac{\mathrm{i}}{\varepsilon}2y\xi}\mathrm{d}y\\
=&-\sum_{j=-\frac{N}{2}}^{\frac{N}{2}-1}\widehat{W}_j(x,t)\frac{1}{(\pi\xi)^d}\int(2\pi)^d\delta(\mu_j+\frac{2y}{\varepsilon})(\delta V)(x,y,t)\mathrm{e}^{-\frac{\mathrm{i}}{\varepsilon}2y\xi}\mathrm{d}y\\
=&-\sum_{j=-\frac{N}{2}}^{\frac{N}{2}-1}\widehat{W}_j(x,t)\frac{1}{(\pi\xi)^d}\int(2\pi)^d(\frac{\varepsilon}{2})^d\delta(y+\frac{\varepsilon}{2}\mu_j)(\delta V)(x,y,t)\mathrm{e}^{-\mathrm{i}y\xi}\mathrm{d}y\\
=&-\sum_{j=-\frac{N}{2}}^{\frac{N}{2}-1}\widehat{W}_j(x,t)(\delta V)(x,-\frac{\varepsilon\mu_j}{2},t)\mathrm{e}^{\mathrm{i}\mu_j(\xi-c)}
=\sum_{j=-\frac{N}{2}}^{\frac{N}{2}-1}\widehat{W}_j(x,t)(\delta V)(x,\frac{\varepsilon\mu_j}{2},t)\mathrm{e}^{\mathrm{i}\mu_j(\xi-c)}	
\end{align*}
where we use the fact:
\begin{align*}
\int \mathrm{e}^{\mathrm{i}kx}\mathrm{d}k=2\pi\delta(x).	
\end{align*}

\section*{Acknowledgement}
The work of Q.Yi was supported by the Training Program for the National Science Foundation of Guangxi Normal University, China (Grant No. 2024PY004). 
L. Xu would like to thank Dr. Hao Wu and Dr. Zhennan Zhou for their generous support and encouragement.

\bibliographystyle{elsarticle-num} 
\bibliography{TSSP_WPFP}

\end{document}